\documentclass{amsart}
\usepackage{amssymb,amsmath}
\usepackage[american,french]{babel}
\usepackage{amsopn}
\usepackage{fancyhdr}

\def\Titre{Exponential Mixing for Stochastic PDEs: The Non-Additive Case. }
\title{\Titre}
\author{Cyril ODASSO
\\
 \\  Ecole Normale Sup\'erieure de Cachan, antenne de Bretagne,\\ Avenue Robert Schuman,
 Campus de Ker Lann, 35170 Bruz (FRANCE). \\ and
\\ IRMAR,  UMR 6625 du CNRS, Campus de Beaulieu,  35042 Rennes cedex (FRANCE)\\}

\newtheorem{Theorem}{Theorem}[section]

\newtheorem{Proposition}[Theorem]{Proposition}
\newtheorem{Lemma}[Theorem]{Lemma}
\newtheorem{Corollary}[Theorem]{Corollary}
\newtheorem{Remark}[Theorem]{Remark}
\newtheorem{Example}[Theorem]{Example}

\makeatletter
\newif\ifmsbmloaded@

\@addtoreset{equation}{section}

\makeatother

\def\R{\mathbb R}
\def\Q{\mathbb Q}
\def\N{\mathbb N}

\def\C{\mathbb C}
\def\E{\mathbb E}
\def\P{\mathbb P}
\def\Pcal{\mathcal{P}}
\def\Hcal{\mathcal{H}}

\def\Dr{\mathcal D}
\def\Br{\mathcal B}

\def\F{\mathcal F}

\newcommand{\BLANC}[1]{   }

\newcommand{\abs}[1]{\left\vert#1\right\vert}
\newcommand{\norm}[1]{\left\Vert#1\right\Vert}
\newcommand{\eps}{\varepsilon}
\newcommand{\sig}{\sigma}
\renewcommand{\i}{\textrm{i}}
\def \Espace{\renewcommand{\arraystretch}{1.7} }
\newcommand{\carre}{    \begin{flushright}
                $\Box$
            \end{flushright}}

\begin{document}

\selectlanguage{american}

\maketitle

\pagestyle{fancy}

\noindent\textbf{Abstract}:
We establish a general criterion which ensures  exponential mixing of
 parabolic Stochastic Partial Differential Equations (SPDE) driven by a
  non additive noise
which is  white in time and smooth in space.
 We apply this criterion on two representative examples: 2D Navier-Stokes (NS) equations
  and Complex Ginzburg-Landau (CGL) equation with a locally Lipschitz noise.
 Due to the possible degeneracy of the noise, Doob theorem  cannot be applied.
Hence a coupling method is used
in the spirit of
\cite{EMS}, \cite{KS3} and  \cite{Matt}.

\noindent Previous results require assumptions on the covariance of the noise
 which might seem restrictive and artificial.
For instance, for NS and  CGL, the covariance operator is supposed
 to be diagonal in the eigenbasis of the Laplacian and not depending on the high modes
  of the solutions.
The method developed in the present paper gets rid of such
assumptions and only requires
  that the range of the covariance operator contains the low modes.

\


\section*{Introduction}

We investigate ergodic properties of parabolic  Stochastic Partial Differential Equations (SPDE)
 driven by a noise which is  white in time and smooth in space.
 Such systems are difficult to handle
with the standard theory because the phase spaces are infinite dimensional.
Moreover the noise is allowed to be degenerate and the conditions required to apply Doob theorem are not always verified
(see \cite{DPZ2} for the theory of ergodicity when Doob Theorem can be applied).

The idea of compensating the degeneracy of the noise on some
subspaces by dissipativity arguments has been introduced in
\cite{KS2000}, \cite{KS00}, and then in  \cite{BKL}, \cite{EMS} .
In the same spirit, we consider systems which
 have only a finite number of unstable directions. In other words, the unstable manifold
is finite-dimensional.  Dissipative SPDEs such
 as the stochastic
 2D Navier--Stokes (NS) and Complex Ginzburg--Landau (CGL) equations have this structure.
 The main requirement on the noise is that it
  is non degenerate in the unstable directions.
  Later, coupling   methods have been introduced to prove exponential convergence to
equilibrium (see \cite{H}, \cite{KS}, \cite{KS2}, \cite{KS3},  \cite{Matt} and \cite{S}).

These articles mainly deal with additive noises. Only in \cite{Matt}, the noise is allowed
to have some dependence on the solution but it has to be of a very special form - see below for
more details. Moreover, the noise is assumed to be diagonal in the eigenbasis of the linear part
of the equation.

In this article, we wish to get rid of these assumptions.
This requires substantial adaptations in the method, for instance
an auxiliary process is introduced.
We develop a general ergodic criterion which ensures exponential mixing
 of the solution provided
the image of the covariance operator of the noise contains the unstable modes.

Roughly speaking, our method allows to treat SPDEs perturbed by a noise of the type $\phi(u)dW$
 where $u$ is the unknown of the equation and $W$ is the driving noise.
 Denoting by $P_N$  the projection onto the unstable modes,
 our main assumption  is that the range of $\phi(u)$ contains the unstable modes
$P_N H$. We think that this is a very natural condition. Note that with these notations,
the above cited articles treat noises of the type $\phi dW$, $\phi$ being constant and diagonal
and with the main  assumption that the range of $\phi$ contains the unstable modes
$P_N H$. In \cite{Matt} (see also \cite{ODASSO1}), the noise has the form
 $\phi(P_N u)dW$ with $(I-P_N)\phi(P_N u)=0$.

 Our method is very general.
 Given a SPDE, it is sufficient to build an  auxiliary process with
 good properties to apply our method
and establish exponential convergence of the solutions to equilibrium.
The technic to build this process depends on the type of SPDE. In fact, we
distinguish three types of  SPDEs. Examples of the two first types are given by NS and CGL.
 The third type of SPDE is more complicated to treat. It includes weakly damped
 but not strongly dissipative SPDEs.
 An example is the weakly damped Non-Linear Schr\"odinger (NLS) equation (see
 \cite{ODASSO2} for the case of an additive noise).
We will study this equation in a forthcoming article.

 The NS equations describe the time evolution of an incompressible fluid.
 It has been widely studied. Most of the articles cited above have been motivated by the application
 to this equation.

Originally introduced to describe a phase transition in superconductivity \cite{GL}, the CGL
 equation also models the propagation of dispersive non-linear waves in various areas of physics such as hydrodynamics
\cite{Newel1}, \cite{Newel2}, optics, plasma physics, chemical reaction \cite{Huber}...
 The CGL equation arises in the same areas of physics as the non-linear Schr\"odinger (NLS) equation.
 In fact, the CGL
equation  is obtained by adding two viscous terms to the NLS equation. The inviscid limits of the deterministic and
stochastic CGL equation to the NLS equation are established in \cite{Bebouche} and \cite{KSlimit}, respectively.

Ergodicity of the stochastic CGL equation is established in \cite{Marc} when the noise is invertible and in
\cite{H} for the one-dimensional cubic case when the noise is  diagonal, does not
 depend on the solution and is smooth in space. In \cite{ODASSO1}, we have established exponential
mixing of CGL driven by a noise which verifies the additional assumptions mentioned above
 under the $L^2$ or the $H^1$--subcritical conditions.

We hope that the method developed here can be combined with other
recent ideas.
 For instance, in  \cite{HM04}, \cite{MattPar}, the case of NS perturbed by a four dimensional noise is treated.
 Hopefully, a four dimensional noise  depending on the unknown could be studied.
 Another topic of interest is to try to prove exponential mixing in the three dimensional case for the
 transition semigroup constructed in  \cite{DebusscheNS3D}. This latter problem will be treated in a forthcoming paper.

Another topic of interest is the study of more general noise. We
will see in forthcoming papers that our Criterion (Theorem
\ref{ThCrt}) could be extended to Levy type  noises.

The remaining of the article is divided into four sections. First
we define some preliminary definitions. In section $2$, we prove
our general criterion (Theorem \ref{ThCrt}) which states that  a
Markov process converges exponentially fast to equilibrium
provided there exists an auxiliary process which verifies some
properties. In section $3$ and $4$ we apply Theorem \ref{ThNS} to
establish
 exponential mixing of the solutions of NS (Theorem \ref{ThNS})
   and CGL (Theorem \ref{ThCGL}).

\section{Preliminaries}
\

\subsection{Cylindrical Wiener process}
\

Let $U, V, K$ be three separable Hilbert spaces. Here
$\mathcal{L}(U;K)$ (resp $\mathcal{L}_2(U;K)$) denotes the space
of bounded (resp Hilbert-Schmidt) linear operators from the
Hilbert space $U$ to $K$. For instance, the inclusion
\begin{equation}\label{Eq0.0}
L^2\left(\left(0,1\right)^d\right)\subset
H^{-m}\left(\left(0,1\right)^d\right),
\end{equation}
 is Hilbert-Schmidt provided $m>\frac{d}{2}$.

The notion of cylindrical Wiener process is a generalization of
the concept of Brownian motion and is used to model noise.

A process is said to be a cylindrical Wiener process of a Hilbert
space $H$ if there exist an orthonormal basis $(e_n)_n$ of $H$, a
family $(W_n)_n$ of independent brownian motions such that
\begin{equation}\label{Eq0.4}
W=\sum_nW_ne_n.
\end{equation}
It is important to notice that $W(t)$ is in $H$ with probability
zero and that the sum converges almost surely (and for any moment)
in $C((0,t);V$) provided the inclusion $H\subset V$  is
Hilbert-Schmidt. Moreover \eqref{Eq0.4} true for a basis $(e_n)_n$
implies \eqref{Eq0.4} for any orthonormal basis of $H$ (See
\cite{DPZ1}).

 \begin{Example}
Let $d\in \N$. It is well-known that the inclusion
$L^2((0,1)^d)\subset H^{-m}((0,1)^d)$ is Hilbert-Schmidt provided
$m>\frac d2$. So, a cylindrical Wiener process $W$ of
$L^2((0,1)^d)$ is a continuous process of $H^{-m}((0,1)^d)$.
\end{Example}

 \begin{Remark}\label{Rqmes}
  In the rest of this paper, when we will
consider a Hilbert space $U$, we will implicitly fix a space $V$
such that $H\subset V$ is Hilbert-Schmidt. Then, all cylindrical
Wiener processes on $U$ will be considered in the same space $V$.
A random variable $X$ living in a space $E$ will be said to depend
measurably on a cylindrical Wiener process $W$ on $U$ if there
exists a map $f: C(0,\infty;V)\to E$ such that $X=f(W)$.
\end{Remark}

\subsection{Topologies on the set of probability measures}
\

\noindent Given a Polish space $E$, the space $Lip_b(E)$ consists
of all the bounded and Lipschitz real valued functions on $E$. Its
norm is given by
$$
\| \varphi \|_L= \abs \varphi _\infty +L_\varphi, \; \varphi\in
Lip_b(E),
$$
where $\abs\cdot _\infty$ is the sup norm and $L_\varphi$ is the
Lipschitz constant of $\varphi$. The space of probability measures
on $E$ is denoted by $\Pcal(E)$. It can be  endowed with the norm
defined by the total variation
 $$
\norm{\mu}_{var}= \sup  \left\{\abs{\mu(\Gamma)}\;|\; \Gamma \in
\Br(E) \right\},
$$
where we denote by $\Br(E)$ the set of the Borelian subsets of
$E$.
 It is well known that $\norm{.}_{var}$ is the dual
 norm of $\abs{.}_\infty$. We can also use the Wasserstein norm
$$
\| \mu\|_*=\sup_{\varphi\in Lip_b(E),\; \norm{ \varphi }_L\le
1}\abs{\int_E \varphi(u)d\mu(u)},
$$
which is the dual norm of $\|\cdot\|_L$.

 \begin{Remark}\label{Rqmes2}
It is important to notice that the set of borelian subset of
$\Pcal(E)$ is the same for both the total variation and the
Wasserstein norm.

Actually, a mapping  $ x\to \P_x$ taking value in $\Pcal(E)$ is
measurable if and only if the mapping $x\to \int_E f(y)d\P_x(y)$
is measurable for a suitable set of map $f: E\to \R$. For
instance, the set of Lipschitz bounded map is such a suitable set.

In section 3 and 4 below, the measurability of the law of the
process $u$ will be a consequence of the fact that
$$
\E\left(\abs{u(t,W,u_0^1)-u(t,W,u_0^2)} \wedge 1\right)\to 0,
$$
provided $u_0^2\to u_0^1$.
\end{Remark}

\subsection{Couplings}
\


 We here recall some results about the coupling. Coupling is the
 basic key of the proof of our criterion (Theorem \ref{ThCrt}
 below). But no known of coupling is required to apply this
 criterion (as we will see in section 3 and 4).

\noindent Let $(\Lambda_1,\Lambda_2)$ be two distributions on a
two space $(E_i,\mathcal{E}_i)_{i=1,2}$. Let
$(\Omega,\mathcal{F},\P)$ be a probability
 space and let $(Z_1,Z_2)$ be a couple random variables $(\Omega,\mathcal{F}) \to (E_i,\mathcal{E}_i)_{i=1,2}$.
 We say that $(Z_1,Z_2)$
is a coupling  of  $(\Lambda_1,\Lambda_2)$  if
$\Lambda_i=\Dr(Z_i)$ for $i=1,2$. We have denoted by $\Dr(Z_i)$
the law of the random variable $Z_i$.

\noindent  Let $\Lambda$, $\Lambda_1$ and $\Lambda_2$ be three
probability measures on a same space $(E,\mathcal{E})$ such that
$\Lambda_1$ and $\Lambda_2$ are absolutely continuous
 with respect to $\Lambda$. We set
$$
d(\Lambda_1 \wedge
\Lambda_2)=(\frac{d\Lambda_1}{d\Lambda}\wedge\frac{d\Lambda_2}{d\Lambda})
d\Lambda.
$$
This definition does not depend on the choice of $\Lambda$ and we
have
$$
\norm{\Lambda_1-\Lambda_2}_{var}= \frac{1}{2} \int_E  \abs{
\frac{d\Lambda_1}{d\Lambda}-\frac{d\Lambda_2}{d\Lambda}}d\Lambda.
$$
Next result is a fundamental result in the coupling methods (See
for instance \cite{Lindvall}).
\begin{Lemma}\label{lem_coupling}
Let $(\Lambda_1,\Lambda_2)$ be two probability measures on a same
$(E,\mathcal{E})$. Then
$$
\norm{\Lambda_1-\Lambda_2}_{var}= \min \P(Z_1\not = Z_2).
$$
The minimum is taken over all couplings $(Z_1,Z_2)$ of
$(\Lambda_1,\Lambda_2)$. There exists a coupling which reaches the
minimum value. It is called a maximal coupling and has the
following property:
$$
\P(Z_1=Z_2, Z_1 \in \Gamma)=(\Lambda_1 \wedge \Lambda_2)(\Gamma)\;
\textrm{ for any } \Gamma \in \mathcal{E}.
$$
\end{Lemma}
\noindent It is interesting to remark that if $\Lambda_1$ is
absolutely continuous with respect to $\Lambda_2$, we have
\begin{equation}\label{e2.34bis}
\norm{\Lambda_1-\Lambda_2}_{var}\leq \frac{1}{2}\sqrt{\int \left (
\frac{d \Lambda_1}{d \Lambda_2} \right )^{2} d\Lambda_2-1}.
\end{equation}
\noindent In order to estimate the bound given in Lemma
\ref{lem_coupling}, we use either
 \eqref{e2.34bis} or the following result which is
lemma D.1 of \cite{Matt} and which is very useful in order to
bound below the probability that a maximal coupling get coupled.
\begin{Lemma}\label{lem_tech_coup_inf}
Let $\Lambda_1$ and $\Lambda_2$ be two equivalent probability
measures on a space $(E,\mathcal{E})$.
  Then for any $p>1$ and any event  $A$ of $E$
 $$
I_p=\int_A \left ( \frac{d \Lambda_1}{d \Lambda_2} \right )^{p}
d\Lambda_1 <\infty \quad \textrm{implies} \quad
\left(\Lambda_1\wedge\Lambda_2\right)(A) \geq
 \left( 1-\frac{1}{p} \right)\left( \frac{\Lambda_1(A)^p}{pI_p} \right)^{\frac{1}{p-1}}.
$$
\end{Lemma}
\noindent Next result is a refinement of Lemma \ref{lem_coupling}
used in \cite{Matt}.
\begin{Proposition}\label{Prop_Matt}
Let $E$ and $F$ be two Polish spaces, $f_0:E\to F$ be a measurable
map  and $(\Lambda_1,\Lambda_2)$ be two probability measures on
$E$. We set
$$
\lambda_i=f_0^* \Lambda_i, \quad i=1,2.
$$
Then there exists a coupling $(V_1,V_2)$ of
$(\Lambda_1,\Lambda_2)$ such that $(f_0(V_1),f_0(V_2))$ is a
maximal coupling of
 $(\lambda_1,\lambda_2)$.
\end{Proposition}
Setting $f_0:(u,v)\to u$, $V_1=(U_1,\widetilde U)$ and
$V_2=(U_2,U_2)$, it follows.
\begin{Corollary}\label{Cor_New}
Let $E$  be a Polish space, $(\Omega,\F,\P)$ be a probability
space and $(U_1,U_2,\widetilde U)$ be three random variables on
$(\Omega,\F,\P)$ taking value in $E$.

 Then there exists a triplet $(u_1,u_2,\widetilde
u)$  such that $(u_2,\widetilde u)$ is a maximal coupling of
$(\Dr(U_2),\Dr(\widetilde U))$ and such that the law of
$(u_1,\widetilde u)$ is $\Dr(U_1,\widetilde U)$.
\end{Corollary}

 \begin{Remark}[Measurability and Markov property]\label{Rqmes3}
\

Let $(u_0^1,u_0^2)\to (U(u_0^1),\widetilde U(u_0^1,u_0^2))$ such
that $(u_0^1,u_0^2)\to \Dr(U(u_0^1),\widetilde U(u_0^1,u_0^2))$ is
measurable (See Remark \ref{Rqmes2}). It is possible to build the
triplet $(u_1,u_2,\widetilde u)$ of Corollary \ref{Cor_New}
associated to $(U(u_0^1),U(u_0^2),\widetilde U(u_0^1,u_0^2))$ such
that $(u_0^1,u_0^2)\to (u_1,u_2,\widetilde u)$ is measurable (See
for instance Remark A.1 of \cite{ODASSO1}).

Now, assume that $u_0^1,u_0^2, u_1,u_2$ live in the same space and
let $((u_0^1,u_0^2)\to(u_1^n,u_2^n,\widetilde
u^n)(u_0^1,u_0^2))_{n}$ be a sequence of independent versions of
the triplet $(u_0^1,u_0^2)\to(u_1,u_2,\widetilde u)(u_0^1,u_0^2)$.
Since it measurably depends on $(u_0^1,u_0^2)$, we can iterate
this sequence, i.e. we set
$$
\Espace\left\{ \begin{array}{rcl}
 (u_1(0),u_2(0))&=&(u_0^1,u_0^2),\\
 (u_1(n+1),u_2(n+1),\widetilde u(n+1))&=& (u_1^{n+1},u_2^{n+1},\widetilde
 u^{n+1})(u_1(n),u_2(n)).
 \end{array}
\right.
$$
It easily follows that $((u_1(n),u_2(n),\widetilde u(n)))_n$ is a
Markov chain.

Another viewpoint of what we have done is the following. We first
have set
$$
(u_1(0),u_2(0))=(u_0^1,u_0^2).
$$
Then, assuming that $(u_1,u_2,\widetilde u)$ is build on
$\{0\dots,n\}$, we have fixed a path and we have build
$(u_1(n+1),u_2(n+1),\widetilde u(n+1))$ as a triplet of
$(U(u_1(n)),U(u_2(n)),\widetilde U(u_1(n),u_2(n)))$. Finally, we
have integrated the probability over the path on $\{0,\dots,n\}$.
This is this viewpoint we will use in the next sections.
\end{Remark}

\section{ A general criterion  }\label{SecCrt}
\

This section is devoted to the statement and the proof of a
general criterion -Theorem \ref{ThCrt}- which ensures exponential
mixing of a Markov process $u$, provided there exists an auxiliary
process $\widetilde u$ which verifies some properties.
 In particular, Theorem \ref{ThNS} below (resp Theorem \ref{ThCGL}) that states that
  the solutions of NS (resp CGL) are exponentially mixing is a Corollary of Theorem
 \ref{ThCrt}.

\subsection{ Statement of the criterion}\label{Sec2.1}
\

 Let $(U,\abs\cdot_U)$  be a Hilbert space and $W$ be a
cylindrical Wiener process on $U$. We are concerned with a
continuous homogenous weak Markov process $u$ taking value in a
Polish space $(H,d_H)$. This Markov process $u$ is assumed to be a
non anticipative measurable map of $W$ (See Remark \ref{Rqmes}).
Since $u$ is a Markov process, it is assumed that its law $\Dr(u)$
is measurably depending of its initial condition $u_0$ (See Remark
\ref{Rqmes2}). We will denote the dependance of $u$ with respect
to $(t,W,u_0)$ as follows
$$
u(t)=u(t,W,u_0).
$$
We denote by $(\mathcal P_t)_{t\in\R^+}$ the Markov transition
semi-group associated to the Markov
 family $\left(u(\cdot,W,u_0)\right)_{u_0\in H}$.

 We first assume  the existence of an auxiliary process
$\widetilde u$ such that $(u,\widetilde u)$ is Markov.

\noindent {\bf A0} {\it There exists a continuous process
$\widetilde u$ taking value in $H$ and that is a non anticipative
measurable map of $W$. Moreover $(u,\widetilde u)$ is a homogenous
weak Markov process and its law $\Dr(u, \widetilde u)$ is
measurably depending of its initial condition $(u_0,\widetilde
u_0)$.}

 We
will denote the dependance of $(u,\widetilde u)$ with respect to
$(t,W,u_0,\widetilde u_0)$ as follows
$$
(u(t),\widetilde u(t))=(u(t,W,u_0),\widetilde u(t,W,u_0,\widetilde
u_0)).
$$

 The next assumptions involve a positive measurable
functional $\Hcal: \,H\to \R^+$, which plays the role of a
Lyapunov
 functional. We assume that there exist $\gamma, C_1, C>0$ and a mapping $h:H^2\to U$
  such that the following hold.

\noindent {\bf A1} {\it There exists  a family
$(C_\alpha')_{\alpha\in (0,\infty)}$, such that for any $u_0\in
H$, any $t\geq 0$, any $\alpha>0$  and any stopping time $\tau\geq
0$,}
$$
\Espace\left\{\begin{array}{rcl}
\E\left(\Hcal(u(t,W,u_0))\right)&\leq&  e^{-\gamma t} \Hcal(u_0)+C_1,\\
\E\left(e^{-\alpha\tau}\Hcal(u(\tau,W,u_0))1_{\tau<\infty}\right)&\leq&
\Hcal(u_0)+C_\alpha'.
\end{array}\right.$$

\noindent {\bf A2} {\it For any $(u_0^1,u_0^2)\in H^2$, for any
couple $(W_1,W_2)$ of cylindrical
 Wiener processes of $U$ and for any
 $t\geq 0$, we have
$$
\P\left(d_H(u_1(t), u_2(t)) \geq C
 e^{-\gamma t}\textrm{ and  } \,\widetilde u=u_2 \textrm{ on } [0,t] \right)
\leq  C e^{ -\gamma t},
$$
where
$$
\Espace \left\{
\begin{array}{rcl}
 u_i(t)&=&u(t,W_i,u_0^i)\quad \textrm{ for }\quad i=1,2,\\
\widetilde u(t)&=&\widetilde u(t,W_1,u_0^1,u_0^2),\\
2C_1&\geq&\Hcal(u_0^1)+\Hcal(u_0^2).
\end{array}
\right.
$$
}

\noindent {\bf A3} {\it For any $(t,u_0^1,u_0^2)\in
(0,\infty)\times H^2$, we have almost surely}
$$
\widetilde u\left(t,W,u_0^1,u_0^2\right)=u\left(t,W+\int_0^\cdot
h\left(u(s,W,u_0^1),\widetilde u(s,W,u_0^1,u_0^2)\right)ds\,,
u_0^2\right).
$$

\noindent {\bf A4} {\it For any couple $(W_1,W_2)$ of cylindrical
Wiener processes of $U$ and for any $(t_0,u_0^1,u_0^2)\in
[0,\infty)\times H^2$
$$
\P\left(\int_{t_0}^\tau  \abs{h(t)}_U^2dt \geq C e^{-\gamma
t_0}\textrm{ and } \,\widetilde u=u_2 \textrm{ on } [s,\tau]
\right) \leq  Ce^{ -\gamma t_0},
$$
 where $(\widetilde u, u_2)$ are defined in {\bf A2}, where $\tau\geq t_0$ is any stopping time and where
$$
\Espace \left\{
\begin{array}{rcl}
h(t)&=& h\left(u(t,W,u_0^1),\widetilde u(t,W,u_0^1,u_0^2)\right),\\
2C_1&\geq&\Hcal(u_0^1)+\Hcal(u_0^2).
\end{array}
\right.
$$}

\noindent {\bf A5} {\it There exists $p_1>0$ such that for any
$(u_0^1,u_0^2)\in H^2$, we have
$$
\P\left(\int_{0}^\infty  \abs{h(t)}_U^2dt \leq C\right) \geq p_1,
$$
 where
$$
\Espace \left\{
\begin{array}{rcl}
h(t)&=& h\left(u(t,W,u_0^1),\widetilde u(t,W,u_0^1,u_0^2)\right),\\
2C_1&\geq&\Hcal(u_0^1)+\Hcal(u_0^2).
\end{array}
\right.
$$}

We now state our criterion.

\begin{Theorem}\label{ThCrt}
Under the above assumptions, there exists a unique stationary
 probability measure $\mu$ of $(\Pcal_t)_{t\in \R^+}$ on $H$.
 Moreover, $\mu$ satisfies
\begin{equation}\label{EqCrta}
\int_{H} \Hcal(u) d\mu(u) <\infty,
\end{equation}
and there exist $C,\gamma'>0$  such that for any $\lambda \in
\Pcal(H)$
\begin{equation}\label{EqCrtb}
\|\Pcal^*_t\lambda-\mu\|_*\leq C e^{-\gamma' t}\left(1+ \int_{H}
\Hcal(u) d\lambda(u) \right).
\end{equation}
\end{Theorem}

\noindent

\noindent Theorem \ref{ThCrt} is proved in sections 2.2, ..., 2.9
hereafter. Let us quickly sketch the proof.

\noindent Assumption {\bf A1} is standard and ensures that
\eqref{EqCrta} holds and that the time of return of the process in
any ball of radius greater than $2C_1$ admits  an exponential
moment.

\noindent Assumption {\bf A2} states that $u(t,W,u_0^1)$ and
$\widetilde u(t,W,u_0^1,u_0^2)$
 become close exponentially fast in probability.

\noindent  Assumption {\bf A3} means that the law $\widetilde
u(\cdot,W,u_0^1,u_0^2)$ is the law of $u(\cdot,W',u_0^2)$ where
$W'$ is a drifted Wiener process. Assumptions {\bf A4} and {\bf
A5}  imply that the Novikov condition holds for a truncation of
$\widetilde u$. So combining  {\bf A3}, {\bf A4} and {\bf A5}, a
Girsanov Transform can be used to build  a couple of Wiener
processes  $(W_1,W_2)$ such that
$$
\widetilde u(\cdot,W_1,u_0^1,u_0^2)= u(\cdot,W_2,u_0^2),
$$
with a positive probability.

 \noindent {\it Conclusion:} Iterating and combining the three properties, we can conclude by
 remarking that it allows to control the probability
of $ u(t,W_1,u_0^1)$ and $ u(t,W_2,u_0^2)$ being very close.
Actually, we wait for entering the ball of radius $2C_1$. Then,
the probability that $(\widetilde u(\cdot-s,W_1,u_1(s),u_2(s)),
u(\cdot-s,W_2,u_2(s)))$ get coupled for any time is bounded below.
If it fails to couple, we wait again for entering  the ball of
radius $2C_1$ and we retry. It follows that there exists a random
time $T_*$ with exponential moment such that
$$
\widetilde u(\cdot-T_*,W_1,u_1(T_*),u_s^2)=
u(\cdot-T_*,W_2,u_2(T_*)).
$$
So applying {\bf A2}, we are able to conclude.

\begin{Remark}
 For the Navier-Stokes and the Complex Ginzburg-Landau equations treated
in section 3 and 4, we use the same functional
$\Hcal=\abs\cdot^2$. More complicated choices may be necessary as
 in
the case of the weakly damped Non Linear Schr\"odinger equation
treated in \cite{ODASSO2}.
\end{Remark}

\begin{Remark}[Finest assumption] It is possible to refine our
criterion as follows. In our proof and especially in section
\ref{Sec2.3}, assumption {\bf A5} is used to prove  the following
irreducibility argument

\

\noindent {\bf A6} For any $r_0>0$, there exist $T,p_0>0$ and a
measurable family of coupling $(u_0^1,u_0^2)\to(u_1,u_2)$ of
$(\Dr(u(\cdot,W,u_0^1),u(\cdot,W,u_0^2)))$ such that
$$
\P\left(d_H(u_1(T),u_2(T))\leq r_0\right)\geq p_0,
$$
provided
$$
2C_1\geq\Hcal(u_0^1)+\Hcal(u_0^2).
$$

 \

Actually {\bf A0},\dots,{\bf A4} and {\bf A6} are sufficient to
prove \eqref{EqCrta} and \eqref{EqCrtb} provided the constant C in
our assumptions verifies suitable conditions when
$d_H(u_0^1,u_0^2)\leq r_0$. For NS and CGL, it is easier to prove
directly {\bf A0}, \dots, {\bf A5}. But, for some equation, {\bf
A5} is not true and {\bf A6} can easily be proved. In that case,
the criterion should be adapted.
\end{Remark}

\begin{Remark}
It we replace the term $e^{-\gamma t}$ by a negative power of $t$
in our assumptions, we obtain a Theorem analogous to Theorem
\ref{ThCrt}, but where convergence is polynomial instead of being
exponential.

\noindent We have seen in \cite{ODASSO2} that for the stochastic
Non--Linear Schr\"odinger (NLS) equation the control of the energy
is polynomial. Hence, we think that the polynomial version of our
criterion
 is the good framework when studying weakly damped SPDE.

\noindent Moreover we will see in Remarks \ref{RqNew1} and
\ref{RqNew2} that
 there exist some variations of Theorem \ref{ThNS} and \ref{ThCGL}
  whose convergence is polynomial. To establish such
results, we need the polynomial version of   Theorem \ref{ThCrt}.
\end{Remark}

\begin{Remark}[Levy type noise]
It is possible to obtain a Theorem analogous to Theorem
\ref{ThCrt} where the noise $W$ is replaced by a Levy process
$N(dtdz)$ (or by a couple $(W,N(dtdz))$). But the proof is more
complicated because the Girsanov transform for measure causes some
problems.

\noindent This generalization is required by some SPDEs that
naturally appear in physical problems such as the stochastic
Non--Linear Schr\"odinger (NLS) driven by Levy-Brown noise.
 These results would be established in forthcoming
papers.
\end{Remark}

The remaining of this section is devoted to the proof of  Theorem
\ref{ThCrt}.

\subsection{Building of a coupling of the solutions}\label{Sec2.2}
\

Let us denote by $(u_0^1,u_0^2)$ two initial conditions in $H$ and
by $T$ a positive real number.

 Applying Corollary \ref{Cor_New},
we build a family of independent measurable map (See Remark
\ref{Rqmes3})
$$
\left((u_0^1,u_0^2)\to (u_1^n,u_2^n,\widetilde
u^n)(u_0^1,u_0^2)\right)_{n\in\N},
$$
such that the law of $(u_1^n,\widetilde u^n)(u_0^1,u_0^2)$ is
 $(\Dr(u(\cdot,W,u_0^1)),\Dr(\widetilde u(\cdot,W,u_0^1,u_0^2)))$
and $(u_2^n,\widetilde u^n)(u_0^1,u_0^2)$ is a maximal coupling of
$(\Dr(u(\cdot,W,u_0^2)),\Dr(\widetilde u(\cdot,W,u_0^1,u_0^2)))$
on $(0,T]$.

 We now build the coupling
$(u_1,u_2)$ of $(\Dr(u(\cdot,W,u_0^1)),\Dr(u(\cdot,W,u_0^2)) )$ on
$[0,nT]$ by induction on $n\in \N$.
 Indeed we first set
$$
u_i(0)=u_0^i,\quad i=1,2.
$$
Assume that $(u_1,u_2,\widetilde u)$ is build on $[0,nT]$, we
extend it on $[0,(n+1)T]$ be setting
$$
(u_1,u_2,\widetilde u)(nT+\cdot)=(u_1^n,u_2^n,\widetilde
u^n)(u_1(nT),u_2(nT)).
$$
It easily follows that $(u_1,u_2)$ is a coupling of
$(\Dr(u(\cdot,W,u_0^1)),\Dr(u(\cdot,W,u_0^2)))$ on $(0,\infty)$
and that the triplet $(u_1,u_2,\widetilde u)$ is homogenous weak
Markov at discrete time $T\N$. Since since the strong Markov
property is equivalent to the weak Markov property when working at
discrete times, it means that for any stopping times $\tau\in
T\N\cup\{\infty\}$
$$
\E_{u_0^1,u_0^2} \left(1_{\tau<\infty}f((u_1,u_2,\widetilde
u))o\,\theta_{\tau}\,\left|\,\F_{\tau}\right.\right)=
1_{\tau<\infty}\E_{u_1(\tau),u_2(\tau)}\left(f(u_1,u_2,\widetilde
u)\right),
$$
where we have denoted by $(\theta_t)_t$ the family of shift
operators
$$
\left(f(u_1,u_2,\widetilde u)\right)o\,\theta_t =
f(u_1(t+\cdot),u_2(t+\cdot),\widetilde u(t+\cdot)).
$$

\subsection{Introduction of $l_0$}
\

\noindent In order to apply {\bf A2}, we  define a family of
integer valued random process $(l_0(k))_{k\in \N\cup\{\infty\}}$
which is particularly convenient when deriving properties of the
triplet
$$
l_0(k)=\min\left\{l \leq k \, | \, P_{l,k}\right\},
$$
 where $\min \phi =\infty$ and
$$
(P_{l,k}) \left \{ \Espace
\begin{array}{l}
\widetilde u=u_2 \;\textrm{ on } (lT,kT),\\
\Hcal(u_1(lT))+\Hcal(u_2(lT))\leq 2C_1.
\end{array}
\right.
$$

The interest of $l_0$ comes from the fact that it allows to apply
{\bf A2} 
    that can be rewritten in the following form
\begin{eqnarray}
\E\left(\abs{u_2(t)-u_1(t)}\wedge 1_{l_0(\infty)\leq l}\right)
&\leq& Ce^{-\gamma (t-lT)},\label{New1.2}
\end{eqnarray}
provided $t\geq lT$.

\subsection{Construction of a useful coupling}
\

In subsections \ref{Sec2.3} and \ref{Sec2.4}, we are interested
with the law of $(u_1,u_2,\widetilde u)$ conditioned by
$l_0(k)=0$. For that purpose we fix a path of $(u_1,u_2,\widetilde
u)$ on $[0,kT]$ such that $l_0(k)=0$ and we build
$(u_1,u_2,\widetilde u)$ on $(kT,(k+1)T]$ as in subsection
\ref{Sec2.2} by applying Corollary \ref{Cor_New}. When the path of
$(u_1,u_2,\widetilde u)$ is fixed on $[0,kT]$, the probability is
denoted by $\Q$.

To understand well the link between $\P$ defined in section 2.2
and $\Q$, notice that
\begin{equation}\label{New1.10}
\Q\left((u_1,u_2,\widetilde u)\in
B\right)=\P\left((u_1,u_2,\widetilde u)\in
B\,\left|\,\F_{kT}\right.\right),
\end{equation}
provided $B$ is a measurable subset of $C(kT,(k+1)T;H)^3$.

It is not easy to work directly with the couple $(\widetilde
u,u_2)$. We prefer working with cylindrical Wiener processes
because of a very useful tool, namely the Girsanov transform.
 That is the reason why we
build the following coupling.

We first set
\begin{equation}\label{New1.4}
\Espace\left\{\begin{array}{rcl}
h(t,W)&=&h(u(t-kT,W,u_1(kT)),\widetilde
u(t-kT,W,u_1(kT),u_2(kT))),\\
\tau(W)&=&\inf\left\{t\in (kT,(k+1)T)\,\left|\,\int_{kT}^t\abs{
h(t,W)}^2dt> 2Ce^{-\gamma kT}\right.\right\}.
\end{array}\right.
 \end{equation}
  Then, applying Corollary \ref{Cor_New} to
$$
\left(W,W,W+\int_{kT}^{\tau(W)\wedge\cdot} h(t,W) dt\right),
$$
 we obtain a couple $(W_1,W_2)$ of cylindrical Wiener processes
such that
$$
\left(W_2,W_1+\int_{kT}^{\tau(W_1)\wedge\cdot} h(t,W_1)dt\right),
$$
is a maximal coupling on $[kT,(k+1)T]$.

Since $l_0(k)=0$, it follows from the definition of $l_0(k+1)$
that
$$
\Q\left(l_0(k+1)=0\right)=\Q\left(\widetilde u=u_2 \right).
$$
Applying the maximal coupling property of the couple $(\widetilde
u,u_2)$, it follows from the fact that $(\widetilde
u(\cdot-kT,W_1,u_1(kT),u_2(kT)),u(\cdot-kT,W_2,u_2(kT)))$ is a
coupling of $(\Dr(\widetilde u),\Dr(u_2))$ on $(kT,(k+1)T)$
$$
\Q\left(l_0(k+1)=0\right)\geq \Q\left(\widetilde
u(\cdot-kT,W_1,u_1(kT),u_2(kT))=u(\cdot-kT,W_2,u_2(kT)) \right),
$$
which yields, by {\bf A3},
\begin{equation}\label{New1.9}
\Q\left(l_0(k+1)=0\right)\geq
\Q\left(W_2=W_1+\int_{kT}^{\tau(W_1)\wedge\cdot} h(t,W_1)dt
\textrm{ and } \tau(W_1)=(k+1)T\right).
\end{equation}
Let us set
$$
\Espace \left\{
\begin{array}{rcl}
A&=&\left\{W\,\left|\,\tau(W)=(k+1)T\right.\right\},\\
\Lambda_1&=& \Dr(W),\\
\Lambda_2&=&\Dr\left(W+\int_{kT}^{\tau(W)\wedge\cdot}
h(t,W)dt\right).
\end{array}\right.
$$
Novikov condition is obviously verified. So, Girsanov Transform
gives
$$
\left(\frac{d\Lambda_
1}{d\Lambda_2}\right)(W)=\exp\left(\int_{kT}^{\tau(W)}
h(t,W)dW(t)-\frac12\int_{kT}^{\tau(W)}\abs{ h(t,W)}^2dt\right),
$$
which yields
\begin{equation}\label{New1.7}
 \int\left(\frac{d\Lambda_
1}{d\Lambda_2}\right)^2d\Lambda_1\leq
\E\left(e^{\int_{kT}^{\tau(W)}\abs{ h(t,W)}^2dt}\right)\leq
e^{2Ce^{-\gamma kT}}.
\end{equation}

\subsection{Probability that $l_0(1)=0$}\label{Sec2.3}
\

 In that subsection, we treat the case $k=0$,
i.e. we assume that $l_0(0)=0$. In that case $\P=\Q$.

Applying the maximal coupling property of the couple
$(W_2,W_1+\int_0^{\tau(W_1)\wedge\cdot} h(t,W_1)dt)$, it follows
from \eqref{New1.9} and from Lemmas \ref{lem_coupling} and
\ref{lem_tech_coup_inf} that
\begin{equation}\label{New1.6}
\P\left(l_0(1)=0\right)\geq  \frac14
\left(\int\left(\frac{d\Lambda_
1}{d\Lambda_2}\right)^2d\Lambda_1\right)^{-1}\Lambda_1(A).
\end{equation}
It follows from {\bf A5} that
$$
\Lambda_1(A)\geq p_1.
$$
Combining \eqref{New1.6} and \eqref{New1.7}, we obtain
\begin{equation}\label{New1.8}
\P\left(l_0(1)=0\right)\geq  \frac{p_1}4 e^{-2C}.
\end{equation}

\subsection{Probability that $l_0(k+1)=0$}\label{Sec2.4}
\

Recall that, for the probability $\Q$, we have fixed a path of
$(u_1,u_2,\widetilde u)$ on $[0,kT]$ such that $l_0(k)=0$ and we
have build $(u_1,u_2,\widetilde u)$ on $[kT,(k+1)T]$ as in
subsection \ref{Sec2.2}.

Notice that \eqref{New1.9} can be rewritten as follows
$$
\Q\left(l_0(k+1)\not=0\right)\leq
\Q\left(W_2\not=W_1+\int_{kT}^{\tau(W_1)\wedge\cdot} h(t,W_1)dt
\textrm{ or } \tau(W_1)<(k+1)T\right),
$$
which yields
\begin{equation}\label{New1.11}
\Espace\begin{array}{r} \Q\left(l_0(k+1)\not=0\right)\leq
 \Q\left(W_2=W_1+\int_{kT}^{\tau(W_1)\wedge\cdot} h(t,W_1)dt
\textrm{ and } \tau(W_1)<(k+1)T\right)\\
+\Q\left(W_2\not=W_1+\int_{kT}^{\tau(W_1)\wedge\cdot}
h(t,W_1)dt\right).
\end{array}
\end{equation}
Notice that, by {\bf A3},
\begin{equation}\label{New1.12}
\Espace\begin{array}{r}
\Q\left(W_2=W_1+\int_{kT}^{\tau(W_1)\wedge\cdot} h(t,W_1)dt
\textrm{ and } \tau(W_1)<(k+1)T\right)\leq\quad\quad\quad\\
\Q\left(u(\cdot-kT,W_2,u_2(kT))=\widetilde
u(\cdot-kT,W_1,u_1(kT),u_2(kT))\textrm{ and }
\tau(W_1)<(k+1)T\right).
\end{array}
\end{equation}
Applying the maximal coupling property of the couple
$(W_2,W_1+\int_0^{\tau(W_1)\wedge\cdot} h(t,W_1)dt)$, it follows
from Lemma \ref{lem_coupling} and from \eqref{e2.34bis},
\eqref{New1.7} that
\begin{equation}\label{New1.13}
\Q\left(W_2\not=W_1+\int_{kT}^{\tau(W_1)\wedge\cdot}
h(t,W_1)dt\right)\leq Ce^{-\gamma kT}.
\end{equation}
Combining \eqref{New1.11}, \eqref{New1.12}, \eqref{New1.13}, we
obtain
\begin{equation}\label{New1.14}
\Espace\begin{array}{l}
 \Q\left(l_0(k+1)\not=0\right)\leq
Ce^{-\gamma
kT}+\\
\Q\left(u(\cdot-kT,W_2,u_2(kT))=\widetilde
u(\cdot-kT,W_1,u_1(kT),u_2(kT))\textrm{ and }
\tau(W_1)<(k+1)T\right).
\end{array}
\end{equation}
Notice that, since we have fixed $(u_1,u_2,\widetilde u)$ on
$[0,kT]$, we cannot apply {\bf A4} to bound the right hand side.

So using \eqref{New1.10} and integrating \eqref{New1.14} by
$(u_1,u_2,\widetilde u)|_{[0,kT]}$ over $l_0(k)=0$, we deduce from
{\bf A4} that
\begin{equation}\label{New1.15}
\P\left(l_0(k+1)\not=0 \textrm{ and } l_0(k)=0\right)\leq
Ce^{-\gamma kT}.
\end{equation}

\subsection{Probability that $l_0(\infty)=0$ and time of failure}
\

We first assume that $l_0(0)=0$.

Notice that, since $l_0(k) =0$ implies $l_0(l)=0$ for any
$\infty\geq k\geq l \geq 0$, it follows that
$$
\P\left(l_0(\infty)\not=0\right)\leq
\sum_{k=0}^\infty\P\left(l_0(k+1)\not=0 \textrm{ and }
l_0(k)=0\right).
$$
Combining \eqref{New1.8} and \eqref{New1.15}, we obtain
$$
\P\left(l_0(\infty)\not=0\right)\leq 1- \frac{p_1}4
e^{-2C}+\frac{e^{-\gamma T}}{1-e^{-\gamma T}}C.
$$
It follows that there exists $T_0$ such that $T\geq T_0$ implies
\begin{equation}\label{New1.16}
\P\left(l_0(\infty)=0\right)\geq p_0=\frac{p_1}{8} e^{-2C}.
\end{equation}
We now fix $T=T_0$.

We denote by $\sig$ the time where $l_0$ stop being zero, i.e.
$$
\sig=\inf\left\{n\in
\N\,\left|\,l_0\left(n\right)>0\right.\right\}.
$$
It follows from \eqref{New1.15} that
$$
\P\left(\sig=k+1\right)\leq Ce^{-\gamma kT},
$$
which yields that, for any $\alpha\in (0,\frac{\gamma T} 6)$,
\begin{equation}\label{Abs_e}
\E \left(e^{3\alpha \sig }1_{\sig<\infty}\right)\leq c.
\end{equation}

\subsection{Moment of $l_0(\infty)$}
\

From now, we stop assuming that $l_0(0)=0$.

Notice that it follows from {\bf A1} that the time of return
$\delta$ in the ball of radius $2C_1$ admits an exponential
moment, i.e., there  exist $\alpha\in (0,\frac{\gamma T} 6)$ and
$c>0$ such that
\begin{equation}\label{Abs_f}
\E \left(e^{\alpha\delta }\right)\leq  c
\sqrt{1+\Hcal(u_0^1)+\Hcal(u_0^2)},
\end{equation}
where
$$
\delta=\min\left\{n\in \N \,|\, \Hcal(u_1(nT))+\Hcal(u_2(nT))\leq
2C_1 \right\}.
$$
For a proof, see for instance (1.56) of \cite{ODASSO1} and apply a
Schwarz inequality.

Now we combine $\delta$ and $\sig$ and we iterate them. We set
$$
\Espace \left\{
\begin{array}{rcllrcll}
\delta_0&=&\delta,&&\\
 \sig_{k+1}&=&\infty  & \textrm{ if } \delta_k=\infty ,&   \sig_{k+1}&=& \sig o\,\theta_{\delta_k}+\delta_k&\textrm{ else},\\
  \delta_{k}&=&\infty  & \textrm{ if } \sig_k=\infty ,&  \delta_k&=& \delta o\,\theta_{\sig_k}+\sig_k&\textrm{
  else}.
\end{array}
\right.
$$
Notice that $\delta_0$ is the first time $n$ of having $l_0(n)=n$.
If $\delta_k<\infty$, then $\sig_{k+1}=\infty$ means that
$l_0(\infty)=\delta_k$. Otherwise, $\sig_{k+1}<\infty$ is the
first time $n>\delta_k$ where $l_0(n)>\delta_k$ and $\delta_{k+1}$
is the first time  $n\geq\sig_{k+1}$ such that $l_0(n)\not =
\infty$. Actually, $l_0(\delta_{k+1})=\delta_{k+1}$.

We set
$$
\rho=\sig+\delta o\,\theta_\sig.
$$
It follows  that
$$
\E\left(e^{\alpha \rho}1_{\rho<\infty}\right)=\E\left(e^{\alpha
\sig} 1_{\sig<\infty}\E\left(e^{\alpha  \delta }o\,\theta_{\sig
T}\,\left|\,\F_{\sig T}\right)\right)\right).
$$
 Applying the Markov property and \eqref{Abs_f},
$$
\E\left(e^{\alpha \rho}1_{\rho<\infty}\right)\leq
c\E\left(e^{\alpha  \sig} 1_{\sig<\infty}\sqrt{1+\Hcal(u_1(\sig
T))+\Hcal(u_2(\sig T))}\right),
$$
 and then the Schwarz inequality,
 $$
\E\left(e^{\alpha \rho}1_{\rho<\infty}\right)\leq
c\sqrt{\E\left(e^{3\alpha \sig}1_{\sig<\infty}\right)}
\sqrt{\E\left(e^{-\alpha \sig}\left(1+\Hcal(u_1(\sig
T))+\Hcal(u_2(\sig T))\right)1_{\sig<\infty}\right)}
 $$
We deduce from {\bf A1} and \eqref{Abs_e} 
 that
there exists  $C_0$ such that
\begin{equation}\label{New1.20}
\E\left(e^{\alpha \rho}1_{\rho<\infty}\right)\leq C_0,
\end{equation}
provided $l_0(0)=0$.

Notice that
$$
\delta_k= \delta_{k-1} + \rho o\, \theta_{\delta_{k-1}},
$$
which yields
$$
 \E\left(e^{\alpha\delta_k}1_{\delta_k<\infty}\right)=
\E\left(e^{\alpha\delta_{k-1}}1_{\delta_{k-1}<\infty}
\E\left(e^{\alpha\rho}1_{\rho<\infty}o\,\theta_{\delta_{k-1}}\,\left|\,\F_{\delta_{k-1}T}\right.\right)\right)
$$
 Applying the Markov property, it follows from
\eqref{New1.20} that
$$
\E\left(e^{\alpha\delta_k}1_{\delta_k<\infty}\right)\leq C_0
\E\left(e^{\alpha\delta_{k-1}}1_{\delta_{k-1}<\infty}\right)
$$
which yields, by \eqref{Abs_f},
\begin{equation}\label{New1.19}
 \E\left(e^{\alpha\delta_k}1_{\delta_k<\infty}\right)\leq cC_0^k\left(1+
 \Hcal(u_0^1) +  \Hcal(u_0^2)  \right)
\end{equation}
We set
$$
k_0=\inf\left\{k\in\N\,\left|\,\sig_{k+1}=\infty\right.\right\}.
$$
It follows from \eqref{New1.16} and from the Markov property that
\begin{equation}\label{New1.17}
\P\left(k_0>k\right)\leq (1-p_0)^k.
\end{equation}
Notice that
$$
l_0(\infty)=\delta_{k_0}.
$$
Let $p\in (1,\infty)$. It follows from H\"older inequality,
\eqref{New1.19} and \eqref{New1.17} that
$$
\Espace
\begin{array}{lcl}
\E\left(e^{\frac{\alpha}{p} l_0(\infty)}\right)&\leq& \sum_k
\E\left(e^{\frac{\alpha}{p}
\delta_k}1_{k_0=k}\right),\\
&\leq& \sum_k \left(\E\left(e^{\alpha
\delta_k}1_{\delta_k<\infty}\right)\right)^{\frac{1}{p}}\P\left(k_0=k\right)^{1-\frac{1}{p}},\\
&\leq& c\left(\sum_k
\left(C_0^{\frac{1}{p}}\left(1-p_0\right)^{1-\frac{1}{p}}\right)^k\right)\left(1+\Hcal(u_0^1)+\Hcal(u_0^2)\right).
\end{array}
$$
Choosing $p$ sufficiently high, we obtain $\gamma'\in
(0,\frac{\gamma}{2T})$ such that
\begin{equation}\label{New1.18}
\E\left(e^{\gamma' l_0(\infty)}\right)\leq C \left(1+
 \Hcal(u_0^1) +  \Hcal(u_0^2)  \right).
 \end{equation}

\subsection{Conclusion}
\

Let $t>0$ and $l\in\N$ such that $t\geq lT$. Notice that
$$
\E\left(\abs{u_2(t)-u_1(t)}\wedge 1\right)\leq
\E\left(\abs{u_2(t)-u_1(t)}\wedge 1_{l_0(\infty)\leq l}\right)+
\P\left(l_0(\infty)>l\right).
$$
Setting $l=\lfloor \frac{t}{2T}\rfloor$, it follows from
\eqref{New1.2} and from \eqref{New1.18}
\begin{equation}\label{New1.1}
\E\left(\abs{u_2(t)-u_1(t)}\wedge 1\right)\leq C e^{-\gamma'
t}\left(1+
 \Hcal(u_0^1) +  \Hcal(u_0^2)  \right).
 \end{equation}

Now we are able to conclude.

 It follows from completeness arguments (See for instance \cite{KS3}) that,
 to establish Theorem \ref{ThCrt}, it is sufficient to prove that
\begin{equation}\label{New1.5}
\abs{\Pcal_tf(u_0^2)-\Pcal_tf(u_0^1)}\leq C\norm f_L e^{-\gamma'
t}\left(1+
 \Hcal(u_0^1) +  \Hcal(u_0^2)  \right),
 \end{equation}
for any $f:H\to \R$ bounded Lipschitz and any $(u_0^1,u_0^2,t)\in
H\times H\times (0,\infty)$.

Since $(u_1,u_2)$ is a coupling of
$(\Dr(u(\cdot,W,u_0^1)),\Dr(u(\cdot,W,u_0^2)) )$, it follows that
$\Pcal_tf(u_0^i)=\E f(u_i(t))$ for $i=1,2$ and $t>0$ which yields
$$
\abs{\Pcal_tf(u_0^2)-\Pcal_tf(u_0^1)}\leq 2\norm f_L
\E\left(\abs{u_2(t)-u_1(t)}\wedge 1\right).
$$
Then we deduce Theorem \ref{ThCrt} from \eqref{New1.1}.

\section{The Two-dimensional Navier-Stokes equations}
\
In this section, we investigate properties of the two-dimensional Navier-Stokes (NS) equations
 with Dirichlet boundary conditions. These equations describe the time evolution of an incompressible fluid and
 are given by
\begin{equation}\label{EqIntroNS}
\left\{
\begin{array}{rcl}
du+\nu(-\Delta)u \, dt+(u,\nabla)u\,dt+\nabla p \, dt &=& \phi(u)dW+fdt,\\
                                                \left(\textrm{div } u\right)(t,x) &=& 0,\;\;\;\;\;\;\;\,\textrm{ for } x\in D,\;\; t>s,\\
                                                u(t,x)&=& 0,\;\;\;\;\;\;\;\, \textrm{ for } x\in \partial D, t>s,\\
                                                u(0,x)&=& u_0(x), \;\textrm{ for } x\in D,
\end{array}
\right.
\end{equation}
where $u(t,x)\in \R^2$ denotes  the velocity field at time $t$ and position $x$,
 $p(t,x)$ denotes the pressure, $\phi(u)dW$ is a random external force field acting on the fluid,
 $\nu>0$ is the viscosity of the fluid and  $D$ is an open bounded domain of $\R^2$ with
regular boundary or $D=(0,1)^2$. Also $f(t,x)dt$ is the
deterministic part of the forcing term. For simplicity in the
redaction, we consider the case $f=0$. The generalization to a
square integrable $f$ is easy.

\noindent In Section \ref{Sec1.1}, we rewrite problem \eqref{EqIntroNS} and we state an ergodic result about NS
 (Theorem \ref{ThNS}).
  To establish it, we apply our criterion (Theorem \ref{ThCrt}).

\subsection{  Notations and Main result }\label{Sec1.1}
\

\noindent We denote by $\abs\cdot$ the norm of $L^2(D;\R^2)$ and by
$\norm\cdot$ the norm of $H^1_0(D;\R^2)$. Let $H$ and $ V $ be the closure of the space of
smooth functions on $D$ with compact support and free divergence
for the norm $\abs \cdot$ and $\norm\cdot$, respectively.

\noindent Let $\Pi$ be the orthogonal projection in $L^2(D;\R^2)$ onto the space $H$. Setting
$$
A = \Pi\left(-\Delta\right), \quad  D(A) =  V  \cap H^2(D;\R^2) \;\textrm{ and }\;
B(u)=\Pi \left((u,\nabla)u\right) ,$$
 we can write problem \eqref{EqIntroNS} in the form
\begin{equation}\label{EqNS}
\left\{
\begin{array}{rcl}
du+\nu A u dt+B(u) dt
 &=& \phi(u)dW,\\
                             u(0)&=& u_0,
\end{array}
\right.
\end{equation}
where $W$ is a cylindrical Wiener process on a Hilbert space $U$.

\noindent  In order to have existence and uniqueness of the solution of \eqref{EqNS},
we make the following assumption

\

\noindent {\bf H0} { \it The function $\phi: \,H\to \mathcal{L}_2(U;H)$ is bounded Lipschitz.
 }

\

\noindent We set
$$
B_0=\sup_{u\in H}\norm{\phi(u)}^2_{\mathcal{L}_2(U;H)}
,\quad
L=L_\phi^2
,
$$
where $L_\phi$ 
is the Lipschitz constant of $\phi$.

\indent There exists a unique $H$--valued solution
$u(\cdot,W,u_0)$ of \eqref{EqNS}. This is a continuous homogenous
weak Markov process and a non anticipative measurable map of $W$.
Moreover the law $\Dr(u(\cdot,W,u_0))$  measurably depends of
$u_0$. Existence can be established by using Galerkin
approximation. To prove uniqueness and measurability with respect
to $u_0$, it is possible to prove the last inequality of Remark
\ref{Rqmes2}.
\begin{Remark}
Using Jenji-Krylov arguments, it is possible to prove that there
exists a measurable mapping $(W,u_0)\to u(\cdot,W,u_0)$ with value
in continuous mapping $(0,\infty)\to H$ such that $u(\cdot,W,u_0)$
is the unique solution of \eqref{EqNS}. However, we do not need it
here.
\end{Remark}
 We denote by $(\mathcal P_t)_{t\in\R^+}$ the Markov
transition semi-group associated to the Markov
 family $\left(u(\cdot,W,u_0)\right)_{u_0\in H}$.

\BLANC{\noindent Rewriting \eqref{EqNS} as follows
$$
u(t)=e^{-A t}u_0-\int_0^te^{-A(t-r)}B(u(r)) dr
+\int_0^te^{-A(t-r)}\phi(u(r))dW(r),
$$
and using a contracting fix-point argument, we obtain that under
{\bf H0}, we have existence and uniqueness of the solution of
\eqref{EqNS} in $H$ when $u_0\in H$.
 Moreover, since measurability property is conserved when taking the limit, there exists a measurable map
 $u$ (See Remark \ref{Rqmes})
 such that $u(\cdot,W,u_0)$ is the unique solution of \eqref{EqNS} and such that the second line of {\bf A0}
 is verified. This result ensures
the strong Markov property of the solutions of \eqref{EqNS}. We
denote by $(\Pcal_t)_{t\in \R^+}$ the Markov transition semi-group
associated to the solutions of \eqref{EqNS}.}

\noindent In our computations, we use the following energy
$$
E_u(t)=\abs{u(t)}^2+\nu \int_{0}^t\norm{u(s)}^2 ds.
$$

\noindent It is well-known  that $\left(A,\Dr(A)\right)$ is
 a self-adjoint operator with discrete spectrum. See \cite{Cons}, \cite{Temam}.
We consider $(e_n)_n$ an eigenbasis of $H$ associated to the increasing sequence $(\mu_n)_n$ of eigenvalues
 of $\left(A,\Dr(A)\right)$. We denote by $P_N$ and $Q_N$
the orthogonal projection in  $H$  onto the space $Sp(e_k )_{1\leq N}$ and onto its complementary,
 respectively.

 Now, we make the assumption  which is used to prove the exponential mixing of
$(\Pcal_t)_{t\in \R^+}$.

\

\noindent {\bf H1} {\it  There exist $N\in \N^*$ and a bounded measurable map $g:H\to \mathcal L(H; U)$ such
that for any $u\in H$ }
$$
\phi(u)g(u)=P_N.
$$

\


\begin{Remark}[Sufficient conditions to satisfy  {\bf H1}]\label{RqH1.1}
\

\noindent A sufficient condition to satisfy 
 {\bf H1} is for instance that $U$ is the orthogonal sum $U_1\oplus U_2$ and
 there exist two measurable  maps \mbox{$\left(\phi_i: \,H\to \mathcal{L}_2(U_i;H)\right)_{i=1,2}$}  such that $\phi_2$ verifies {\bf H1}
 and
$$
\phi(u)W=\phi_1(u)W_1+\phi_2(u) W_2,
$$
 for any  $u\in H$ and any $W=(W_1,W_2)\in U$. Moreover, if $\phi_2$ is a constant map, we can omit the
orthogonality condition on $U_1$ and $U_2$. 

\noindent A very interesting consequence is the case of the sum of
a multiplicative noise and an additive noise which covers  the low
modes. Namely, for any
   measurable  map $f:\R^2\to \R$  and $ (\phi_1,\phi_2)\in\mathcal{L}_2(U_1;H^2(D))\times\mathcal{L}_2(U_2;H)$,
 one define $\phi$ by
$$
\phi(u)W=f(u) (\phi_1 W_1)+\phi_2 W_2,
$$
 for any  $u\in H$ and any $W=(W_1,W_2)\in U$.
The operator $\phi$ verifies {\bf H1} provided
$$
P_N H\subset \textrm{Im }\phi_2.
$$

\noindent  Another sufficient condition is that  $U=H$ and  there
exists an invertible operator $\psi$ on the low modes and a
constant $\eps$ such that for any $u\in H$
$$
Q_N\phi(u)P_N=0\quad\textrm{ and }\quad
\abs{P_N\phi(u)P_N-\psi}_{\mathcal L(P_N
H)}<\eps<\abs{\psi^{-1}}^{-1}_{\mathcal L(P_N H)}.
$$
Thus, our result holds when the covariance operator is a small
perturbation of a constant in the low modes.

\noindent We think that these examples might be physically
relevant.
\end{Remark}
\begin{Remark}\label{RemH1}
The existence of a map $\widetilde g(u)$ such that $
\phi(u)\widetilde g(u)=P_N, $ is equivalent to the following
property
$$
P_N H\subset Im(\phi(u)).
$$
\noindent Hence  {\bf H1}  can be seen as a non degeneracy
condition on the low modes in the spirit of
  \cite{BKL}, \cite{ODASSO2}, \cite{EMS}, \cite{H}, \cite{KS00}, \cite{KS},
 \cite{KS2}, \cite{KS3},  \cite{Matt}, \cite{ODASSO1} and \cite{S}.
The lack of surjectivity of $\phi(u)$ on the high modes is
counterbalanced by the dissipativity of \eqref{EqNS}.

\noindent Moreover, if there exists such a map $\widetilde g$,
 then there exists a measurable map $g$ such that
$\phi(u) g(u)=P_N$ and $\abs{g(u)}_{\mathcal L(H;U)}\leq
\abs{\widetilde g(u)}_{\mathcal L(H;U)}$. This mapping $g$ is
constructed by similar ideas as in the construction of the pseudo
inverse.
 Hence the assumption of measurability of $g$ in {\bf H1} is superfluous.
\end{Remark}


\noindent The aim of this section is to prove the following
result.
\begin{Theorem}\label{ThNS}
Assume that {\bf H0} holds. There exists $N_0(B_0,\nu,D,L)$ such that, if  {\bf H1} holds with $N\geq N_0$, then
there exists a unique stationary
 probability measure $\mu$ of $(\Pcal_t)_{t\in \R^+}$ on $H$.
 Moreover, $\mu$ satisfies
\begin{equation}\label{EqThNSa}
\int_{H} \abs{u}^{2} d\mu(u) < \infty,
\end{equation}
and there exist $C,\gamma>0$, such that for any $\lambda \in \Pcal(H)$
\begin{equation}\label{EqThNSb}
\|\Pcal^*_t\lambda-\mu\|_*\leq C e^{-\gamma t}\left(1+ \int_{H} \abs{u}^2 d\lambda(u) \right).
\end{equation}
\end{Theorem}
\noindent The remaining subsections of this section are devoted to
prove Theorem \ref{ThNS} by applying Theorem \ref{ThCrt}.

\noindent The proof of our result is based on coupling arguments. These arguments have initially been used
in the context of dissipative SPDEs in  \cite{KS3}, \cite{Matt}.
For a better understanding of this kind of method, see section 1 of \cite{Matt} and the two first subsections
of \cite{ODASSO1}.
 There the coupling method is explained on two examples which are simpler but contain all the difficulties.

\noindent The method used  in  \cite{KS3} and \cite{Matt} requires
 the three following assumptions for a $N$ sufficiently high.

\noindent The first assumption is a structure condition on $\phi$.  It is a slight generalization of the usual
 assumption that $\phi(u)$
is diagonal in the basis $(e_n)_n$.

\noindent {\bf Ha} {\it The Hilbert space $U$ is $H$. Moreover, for any $u\in H$, we have }
$$
P_N\phi(u)Q_N=0,\; Q_N\phi(u)P_N=0.
$$

\noindent The second assumption means that $\phi$ only depends on its low modes.

\noindent {\bf Hb} {\it For any $u\in H$}
$$
\phi(u)=\phi(P_N u).
$$

\noindent The third assumption is {\bf H1}. Under {\bf Ha}, it could be written in the form.

\noindent {\bf Hc} {\it The linear map $P_N\phi(u)P_N$ is invertible on $P_N H$. Moreover }
$$
\sup_{u
\in H}\abs{\left(P_N\phi(u)P_N\right)^{-1}}_{\mathcal L\left( P_N H;P_N H \right)}<\infty.
$$

\noindent In these papers, the proof is divided in two steps.

\noindent {\it Step 1:} Starting from initial data $(u_0^1,u_0^2)$
in a ball of radius $R_0$,
 Girsanov transform can be used to
show that there exists a couple of Wiener processes $(W_1,W_2)$ such that
$$P_N u(T,W_1,u_0^1)=P_N u(T,W_2, u_0^2)$$
with positive probability.

\noindent  {\it Step 2:} Another ingredient, the so-called Foias-Prodi estimate, is used.
 It is based on the observation that if two solutions are such that
$$(P_N u(\cdot,W_1,u_0^1),Q_N W_1)=(P_N u(\cdot,W_2,u_0^2), Q_N W_2),$$
 for a long time then they become very close. Girsanov transform can again  be used to
show that, if we start from initial conditions which have the same low modes, then
  there exists a couple of Wiener processes  $(W_1,W_2)$ such that
the low modes of the solutions remain equal for all times with a
positive probability.

\noindent {\it Conclusion:}  Since the time of entering a ball of radius $R_0$ admits an exponential moment,
 we are able to combine and to iterate the two steps and then to conclude.

\bigskip

In the above mentioned articles, {\bf Ha} and {\bf Hb} are essential when using the Girsanov transform.
 They are also necessary  to prove that the Foais-Prodi estimate holds pathwise. This latter point is not
 important. We have shown in \cite{ODASSO2} how to use a Foias-Prodi estimate in expectation and it is
 not difficult to see that such an estimate holds in our case.

 \noindent On the contrary, it seems very difficult to use Girsanov Transform without  {\bf Ha} and {\bf Hb}.
  We think that these are artificial whereas  {\bf H1} is very natural.

\noindent The idea used in this article is to separate the use of
the Foias-Prodi estimate and the
 Girsanov Transform by introducing an auxiliary process $\widetilde u(\cdot,W,u_0^1,u_0^2)$
 (In some some sense, this idea could be related to the binding processes of \cite{H}).
 Our criterion (Theorem \ref{ThCrt}) essentially requires that the auxiliary process verifies two properties.

\noindent {\it First fundamental property:}  The first property is a variation of the Foias-Prodi estimate.
   It states that $u(t,W,u_0^1)$ and $\widetilde u(t,W,u_0^1,u_0^2)$
 become close exponentially fast in probability.

\noindent {\it Second fundamental property:}   The second property
is {\bf A3} and
 states that there exists $h$  such that
$$
u(t,W,u_0^1,u_0^2)=u\left(t,W+\int_0^\cdot h(
u(s,W,u_0^1),\widetilde u(s,W,u_0^1,u_0^2))ds\,, u_0^2\right).
$$
 Hence, taking into account {\bf H1}, it follows from coupling arguments and from a Girsanov Transform that there
 exists a couple of Wiener processes  $(W_1,W_2)$ such that
$$
\widetilde u(\cdot,W_1,u_0^1,u_0^2)= u(\cdot,W_2,u_0^2),
$$
with a positive probability.

 \noindent {\it Conclusion:} Iterating and combining the two properties, we can conclude by
 remarking that it allows to control the probability
of $ u(t,W_1,u_0^1)$ and $ u(t,W_2,u_0^2)$ being very close.

\subsection{Energy estimate and Lyapunov structure}\label{SecNSnrj}
\

Setting $\Hcal=\abs\cdot^2$, the following result states that {\bf
A1} is true.
\begin{Lemma}[The Lyapunov structure]\label{lemNSLyap}
Assume {\bf H0}. There exist $C_1$ and a family
$(C_\alpha')_\alpha$ only depending on $\nu, B_0$ and $D$ such
that
$$
\Espace\left\{\begin{array}{rcl}
\E\left(\abs{u(t,W,u_0)}^{2}\right)&\leq&  e^{-\nu \mu_1  t} \abs{u_0}^{2q}+C_1,\\
\E\left(e^{-\alpha
\tau}\abs{u(\tau,W,u_0)}^{2}1_{\tau<\infty}\right)&\leq&
\abs{u_0}^{2}+C_\alpha',
\end{array}\right.$$
for any $t\geq 0$, any $\alpha>0$ and any stopping time $\tau$.
\end{Lemma}

{\it Proof of Lemma  \ref{lemNSLyap}}

\noindent We set $u=u(\cdot,W,u_0)$. We apply Ito Formula to
$\abs{u}^2$
$$
d\abs{u}^2+2\nu\norm{u}^2dt+2\left(u,B(u)\right)dt
=2\left(u,\phi(u)dW\right)+ \norm{\phi(u)}^2_{\mathcal
L_2(U;H)}dt.
$$
Recall that
$$
\left(u,B(u)\right)=0.
$$
Hence, taking into account {\bf H0}, we deduce
\begin{equation}\label{EqA0}
d\abs{u}^2+2\nu\norm{u}^2dt=2\left(u,\phi(u)dW\right)+B_0dt.
\end{equation}
Notice that $\norm\cdot^2\geq \mu_1\abs\cdot^2$. Hence,
integrating and taking the expectation, we establish the first
inequality of Lemma \ref{lemNSLyap}.

Let $\alpha>0$.  Ito Formula of $e^{-\alpha t}\abs{u(t)}^2$ gives
$$
de^{-\alpha s} \abs{u}^2 +e^{-\alpha s}
\left(2\nu\norm{u}^2+\alpha\abs{u}^2\right)dt=2e^{-\alpha
s}\left(u,\phi(u)dW\right)+B_0e^{-\alpha s}dt.
$$
Let $\tau$ be a stopping time and $n\in \N$. Taking the
expectation, it follows
$$
\E\left(e^{-\alpha \tau\wedge n}\abs{u(\tau\wedge n)}^2\right)\leq
\abs{u_0}^2+C_\alpha'.
$$
Let $n\to \infty$, the second inequality of Lemma \ref{lemNSLyap}.
 \carre

The next result will be useful in the following.
 \begin{Lemma}[Exponential estimate for the growth of the solution]\label{lemNSnrj}
Assume that {\bf H0} holds. There exists $\gamma_0>0$ only
depending on $\nu, B_0$ and $D$ such that
$$
\E\left(\exp\left(\gamma_0\sup_{t\geq
0}\left(E_{u(\cdot,W,u_0)}(t)-B_0 t\right)\right)\right) \leq 2
e^{\gamma_0 \abs{u_0}^2},
$$
for any $u_0\in H$.
\end{Lemma}

{\it Proof of Lemma  \ref{lemNSnrj}}

\noindent We set for any $\gamma>0$
$$
M(t)=2\int_0^t\left(u(r),\phi(u(r))dW(r)\right), \quad \mathcal
M_\gamma(t)=M(t)-\frac{\gamma}{2}\left<M\right>(t).
$$
Remarking that
$$
d\left<M\right>=4\abs{\phi(u)^*u}^2dt\leq 4cB_0\norm{u}^2dt,
$$
and setting $\gamma_1=\frac{\nu}{2cB_0}$, we deduce from
\eqref{EqA0} that
\begin{equation}\label{EqA2}
E_u(t)\leq \mathcal M_{\gamma_1}(t)+\abs{u_0}^2+ B_0 t.
\end{equation}
Notice that $e^{\gamma\mathcal M_\gamma}$ is a positive
supermartingale whose value is $1$ at time $0$. We deduce from
maximal supermartingale inequality that
\begin{equation}\label{EqA2bis}
\P\left(\sup_{t\geq 0}\mathcal M_{\gamma_1}(t)\geq \rho\right)
\leq\P\left(\sup_{t\geq 0}e^{{\gamma_1}\mathcal
M_{\gamma_1}(t)}\geq e^{\gamma_1\rho} \right) \leq
e^{-\gamma_1\rho}.
\end{equation}
We set $\gamma_0=\frac{\gamma_1}{2}$. Notice that
$$
\E\left(e^{{\gamma_0}\sup\mathcal M_{\gamma_1}}
\right)=1+\gamma_0\int_0^\infty e^{\gamma_0 x}\P\left(\sup\mathcal
M_{\gamma_1}\geq x\right)dx,
$$
which yields, by \eqref{EqA2bis},
\begin{equation}\label{EqA2.1}
\E\left(e^{\gamma_0\sup \mathcal M_{\gamma_1}}\right)\leq 2.
\end{equation}
Combining \eqref{EqA2} and \eqref{EqA2.1}, it follows
$$
\E\left(\exp\left(\gamma_0\sup_{t\geq 0}
\left(E_{u(\cdot,W,u_0)}(t)-B_0 t\right)\right)\right) \leq
2e^{\gamma_0 \abs{u_0}^2}.
$$
This ends the proof of Lemma  \ref{lemNSnrj}. \carre

\subsection{  Construction of the auxiliary process }\label{Sec1.2}
\

\noindent Now, we build the auxiliary process. We set
$$
F(u)=\nu A(u)+B(u).
$$
Taking into account {\bf H1}, we remark that {\bf A3} is a
consequence of
\begin{equation}\label{EqAux2bis}
\left\{
\begin{array}{rcl}
d\widetilde u+F(\widetilde u) dt + P_N \delta(u(t,W,u_0),\widetilde u)dt&=& \phi(\widetilde u)dW,\\
\widetilde u(0,W,u_0,\widetilde u_0)&=& \widetilde u_0,
\end{array}
\right.
\end{equation}
with
\begin{equation}\label{EqAux2}
h(u,\widetilde u)=-g(\widetilde u)\delta(u,\widetilde u).
\end{equation}
Since we want that $\widetilde u(t,W,u_0^1,u_0^2)$ and $
u(t,W,u_0^1)$ become
 very close in probability,
it is natural to build $\widetilde u$ such that \eqref{EqAux2bis}
and \eqref{EqAux2} hold with
\begin{equation}\label{EqAux3}
\delta(u,\widetilde u)=K   P_N(\widetilde u-u).
\end{equation}
Hence we consider the following equation
\begin{equation}\label{EqNSbis}
\left\{
\begin{array}{rcl}
d\widetilde u+F(\widetilde u) dt +K   P_N(\widetilde u-u(t,W,u_0^1))dt&=& \phi(\widetilde u)dW,\\
                             \widetilde u(0)&=& \widetilde u_0,
\end{array}
\right.
\end{equation}
where $K>0$  will be chosen later and $N$ is the integer used in
{\bf H1}.

\noindent As for \eqref{EqNS}, we deduce from {\bf H0}  that there
exists a unique $H$--valued solution to \eqref{EqNSbis} when
$(u_0,\widetilde u_0)\in  H^2$. Moreover, $\widetilde u$  is a non
anticipative measurable map of $W$ (See Remark \ref{Rqmes}). It
follows from the uniqueness that $(u,\widetilde u)$ is a
homogenous weak Markov process. Moreover, its law $\Dr(u,
\widetilde u)$ is measurably depending of its initial condition
$(u_0,\widetilde u_0)$ (To prove it, one can establish the last
inequality of Remark \ref{Rqmes2}: Arguments are similar but
simpler than those of the proof of the next proposition).

 We
will denote the dependance of $(u,\widetilde u)$ with respect to
$(t,W,u_0,\widetilde u_0)$ as follows
$$
(u(t),\widetilde u(t))=(u(t,W,u_0),\widetilde u(t,W,u_0,\widetilde
u_0)),
$$
and we deduce {\bf A0}.

\noindent Taking into account {\bf H1} and the uniqueness of the
solution of \eqref{EqNS} under {\bf H0},
 we deduce {\bf A3} with \eqref{EqAux2} and \eqref{EqAux3}.

\noindent Now, it remains to prove {\bf A2}, {\bf A4} and {\bf
A5}, i.e. that $\widetilde u(t,W,u_0,\widetilde u_0)$ and $
u(t,W,u_0)$ become
 close exponentially fast.  In   \cite{KS3}, \cite{Matt} and \cite{ODASSO1}, a pathwise Foias-Prodi is used.
Here,  $u-\widetilde u$ does not seem to tend to $0$, pathwise.
That is the reason why we adapt an idea we have already used in
\cite{ODASSO2} to prove polynomial mixing for the weakly damped
Non-linear Schr\"odinger (NLS) equations. Since it seems that
there is no pathwise Foias-Prodi estimate for NLS,   a
Foias-Prodi estimate in expectation
 was used.  Using analogous technics, we have the following result.
\begin{Proposition}\label{PropNSFP}
Assume that {\bf H0} holds. There exist $\gamma>0$, $\eps\in
(0,1]$, $K_0>0$ and $N_0$ depending only on $B_0$, $\nu$ and $D$
such that
 for any $K\geq K_0$ and $N\geq N_0$ and any $(t,u_0,\widetilde u_0)\in (0,\infty)\times H^2$
$$
\E\left(\left(e^{ t}\abs{r(t)}^2+\int_0^t e^{
s}\norm{r(s)}^2ds\right)^\eps\right) \leq
2\abs{r(0)}^{2\eps}e^{\gamma \abs{u_0}^2},
$$
where
$$
r=\widetilde u(\cdot,W,u_0,\widetilde u_0)- u(\cdot,W,u_0).
$$
\end{Proposition}
Notice that, by Chebyshev inequality, this result obviously
implies {\bf A2}, {\bf A4} and {\bf A5}.

{\it Proof of Proposition \ref{PropNSFP}}

\noindent For any function $f$, we denote by $\delta f(u)$ the
value $f(\widetilde u)-f(u_1)$. Taking the difference between
\eqref{EqNS} and \eqref{EqNSbis}, we obtain
$$
dr + \nu Ar dt +KP_N  r dt +\delta B(u)dt
=\delta \phi(u)dW.
$$
Hence, applying Ito Formula to $\abs{r}^2$, we have
\begin{equation}\label{EqA1}
d\abs{r}^2 + 2\left(\nu\norm{r}^2 +K\abs{P_N r}^2\right)dt=
2\left(r,\delta \phi(u)dW\right)+ I(t)dt,
\end{equation}
where
$$
I(t)=I_1(t)+I_2(t)
,\quad 
I_1(t)=-2\left(r,\delta B(u)\right),\quad
I_2(t)=\norm{\delta\phi(u)}^2_{\mathcal L_2(U;H)}.
$$
Remarking that
$$
\delta B(u)=\pi\left( ( \widetilde u,\nabla)r+(r,\nabla) u
\right),\quad (r,(\widetilde u,\nabla) r)=0,
$$
we deduce from a Schwartz inequality that
$$
I_1(t)\leq c \abs{r^2}\norm{u}=c \abs{r}^2_4\norm{u}.
$$
It follows from Sobolev Embedding
$H^{\frac{1}{2}}\left(D,\R^2\right)\subset
L^{4}\left(D,\R^2\right)$ and interpolatory inequality that
$$
I_1(t)\leq c \norm{r}^2_\frac{1}{2}\norm{u_1}\leq c
\norm{r}\abs{r}\norm{u}.
$$
We infer from an arithmetico-geometric inequality that
\begin{equation}\label{EqA1ter}
I_1(t)\leq \nu\norm{r}^2+c\abs{r}^2\norm{u}^2.
\end{equation}
Applying {\bf H0}, we obtain
\begin{equation}\label{EqA1bis}
I_2(t)\leq L\abs{r}^2.
\end{equation}
Remarking that
$$
\left(K\wedge (\nu\mu_{N+1})\right)\abs{r}^2\leq\nu\norm{Q_N
r}^2+K\abs{P_N r}^2\leq\nu\norm{r}^2+K\abs{P_N r}^2,
$$
we deduce from \eqref{EqA1},
 \eqref{EqA1ter} and \eqref{EqA1bis} that there exists $\Lambda$ such that
$$
d\abs{r}^2 + \left(\left(K\wedge
(\nu\mu_{N+1})-L\right)-\Lambda\norm{u_1}^2\right)\abs{r}^2dt\leq
2\left(r,\delta \phi(u)dW\right).
$$
Integrating this formula and taking the expectation, it follows
$$
\E\left(e^tG(t)^{-1}\abs{r(t)}^2+\int_0^t
e^sG(s)^{-1}\norm{r(s)}^2ds\right) \leq \abs{r(0)}^{2},
$$
where
$$
G(t)=e^{ -\left(K\wedge
(\nu\mu_{N+1})-L-1\right)t+\Lambda\int_0^t\norm{u(r)}^2dr}.
$$
It follows from H\"older inequality that
$$
\Espace\begin{array}{r}
 \E\left(\left(e^{ t}\abs{r(t)}^2+\int_0^t
e^{ s}\norm{u(s)}^2ds\right)^\eps\right) \leq \sqrt{\E\left(\sup
G^{2\eps} \right)}\times\quad\quad\quad\\
\left(\E\left(e^tG(t)^{-1}\abs{r(t)}^2+\int_0^t
e^sG(s)^{-1}\norm{r(s)}^2ds\right)\right)^{\eps}.
\end{array}
$$
Choosing $N$, $K$ sufficiently high and $\eps>0$ sufficiently
small, it follows from Lemma \ref{lemNSnrj} that
$$
\E\left(\sup G^{2\eps}\right)\leq 2e^{\gamma_0 \abs{u_0}^2}
$$
which yields Proposition \ref{PropNSFP}. \carre

\subsection{Conclusion and Remarks}\label{sec_Conclusion_NS}
 \

We have proved that assumptions of Theorem \ref{ThCrt} are
verified. So Theorem \ref{ThNS} follows.

Actually, {\bf A0} is a consequence of the well-posedness of
equations.

We set $\Hcal=\abs\cdot^2$.  As shown in Lemma \ref{lemNSLyap}, we
have {\bf A1}.

Assumption {\bf A3}  has been proved at the beginning of
 section 3.3.

We deduce {\bf A2}, {\bf A4} and {\bf A5} directly from
Proposition \ref{PropNSFP}

  Since {\bf A0},\dots,
{\bf A5} are established, we can apply Theorem \ref{ThCrt} which
yields Theorem \ref{ThNS}.

\begin{Remark}\label{RqNew1}
Theorem \ref{ThNS} could be improved and the assumptions could be
weakened. We chose to restrict
 to this statement for clarity and readability.
  For instance, it is possible to replace {\bf H1} by

\noindent {\bf H1'} {\it  There exist $n\in\N^*$, a measurable map
$g:H\to \mathcal L(H; U)$ and two constants $\sig, C$
 such that for any $u\in H$ }
$$
\phi(u)g(u)=P_N,\quad \abs{g(u)}_{ \mathcal L(H; U)}\leq
C\exp\left(\sig \abs{u}^2\right).
$$
In this case $N_0$ depends on $\sig$. Moreover, it is easy to
strengthen \eqref{EqThNSa} into
$$
\int_{H} \exp\left(\sig_1(B_0,\nu,D) \abs{u}^2\right) d\mu(u) <
\infty.
$$

\noindent In {\bf H0}, the boundedness of $\phi$ could be replaced
by $\abs{\phi(u)}\leq C(1+\abs{u}^\gamma)$ with $\gamma<1$. In
this case, the rate of convergence becomes greater than any power
of time instead of being exponential in time. Moreover for any $p$
there exists $c_p$ such that if there exists $C$ such that
$\abs{\phi(u)}\leq C+c_p\abs{u}$, then the rate of convergence is
greater than $(1+t)^{-p}$ instead of being exponential.

\noindent Assume now that $D$ is the two-dimensional torus. We
replace the Dirichlet boundary condition by periodic condition and
we assume that {\bf H0} and {\bf H1'} hold for $N$ sufficiently
high and
$$
B_1=\sup_{u\in H}\norm{\phi(u)}^2_{\mathcal{L}_2(U; V )}<\infty,
$$
\noindent then we can strengthen \eqref{EqThNSa} by
$$
\int_{H} \exp\left(\sig_2(B_1,\nu,D) \norm{u}^2\right) d\mu(u) <
\infty,
$$
and  in \eqref{EqThNSb} we can replace $\norm\cdot_*$, the
Wassertein norm in $H$, by the Wassertstein norm in $H^s(D;\R^2)$
for any $s<1$. Moreover, if $\phi: H\to \mathcal L_2(U; V )$ is
bounded Lipschitz,
 then in \eqref{EqThNSb} we can replace $\norm\cdot_*$ the Wassertein norm in $H$ by the Wassertstein norm in
$ V $.
\end{Remark}

\begin{Remark}
Theorem \ref{ThNS} is of the same type as the results obtained in
\cite{EMS}, \cite{KS3}, \cite{Matt} and \cite{ODASSO1}
 where  SPDEs with additive noise only depending on the low modes are studied.
 Since the decrease of $\Pcal^*_t\lambda-\mu$ is measured in Wasserstein norm, we know that $\Pcal_t \phi$
 converges to its average with respect to  $\mu$ for any Lipschitz function $\phi$.
 In fact, in the above mentioned articles, it is also true  if $\phi$ is only Lipschitz with respect to the high modes.
  We do not know if this holds in our situation.
\end{Remark}

\begin{Remark}\label{Rq_PropNS}
Proof of Proposition \ref{PropNSFP} displays a nice and well known
property of Navier-Stokes equations. Indeed, we see that the
difference between $u$ and $\widetilde u$ can be estimated using
only the energy of $u$ and not  the energy of $\widetilde u$ which
is much more difficult to estimate. No control on the
 probability of the linear growth of the energy of the auxiliary process is required.
 This property holds also for the one-dimensional Burgers equation and
 for the Complex Ginzburg-Landau equation with a globally Lipschitz noise in the subcritical case. However
  it does not hold in general and in this case the construction of the auxiliary process is more involved.
   For instance, in the case of the Non-Linear Schr\"odinger equations, it is not possible to prove a result
   similar to Proposition \ref{PropNSFP}. It is important to show that our method can also be used for these
   equations. The  Non-Linear Schr\"odinger equations will be treated in a forthcoming article and we consider
   the Complex Ginzburg-Landau with a locally Lipschitz noise in section 4.
\end{Remark}


\section{The Complex Ginzburg--Landau equation with a locally Lipschitz noise}
\

\noindent The aim of this section is to apply our method to the
stochastic CGL equation
 with Dirichlet boundary conditions and with a locally Lipschitz noise.

\noindent Let us recall that it has the form
\begin{equation}\label{EqIntroCGL}
\left\{
\begin{array}{rcl}
du+(\eps+\i)(-\Delta)u \,dt+(\eta+\lambda\i)\abs{u}^{2\sig} u \,dt &=& \phi(u)dW+f dt,\\
        u(t,x)&=& 0, \;\;\;
        \textrm{for } x\in \partial D, t>0,\\
        u(0,x)&=& u_0(x),  \textrm{ for } x\in D,
\end{array}
\right.
\end{equation}
where $D$ is an open bounded domain of $\R^d$ with regular boundary or $D=(0,1)^d$,
 where $\eps >0$, $\eta >0$, $\lambda\in\{-1,1\}$
and where we impose the $L^2$--subcritical condition $\sig  d< 2$.
For simplicity in the redaction, we consider the case $f=0$, where
 $f$ is the deterministic part of the forcing term $\phi(u)dW+fdt$.
 The generalization to a square integrable $f$
 is easy.

\noindent Ergodicity for the stochastic CGL equation is established in \cite{Marc} when the noise is invertible and in
\cite{H} for the one-dimensional cubic case when the noise is  diagonal, does not
 depend on the solution and is smooth in space. Then, in \cite{ODASSO1}, we have established exponential
mixing of CGL driven by a noise which verifies {\bf Ha}, {\bf Hb} and {\bf Hc}
 under the $L^2$ or the $H^1$--subcritical conditions.

\noindent As explained  in Remark \ref{Rq_PropNS}, technics of
section 3 can easily be applied to the stochastic CGL equation
with a globally Lipschitz noise. It gives exponential mixing in
$L^2$ under the $L^2$--subcritical condition $\sig d<2$. Moreover,
one can obtain
 the exponential mixing in  $H^1$ under
the $H^1$--subcritical condition $(d-2)\sig <2$ when $\lambda=1$.
Using a polynomial version of our criterion,
 one can obtain polynomial mixing in  $H^1$ under the $L^2$--subcritical condition $\sig d<2$ when $\lambda=-1$.

As explained in Remark \ref{Rq_PropNS}, it seems that such technics can not always be applied when
there is no analogous property to Proposition \ref{PropNSFP}.
For instance, the case of the stochastic non-linear Schr\"odinger equation requires more sophisticated tools
 and will be treated in a forthcoming paper.
We study the CGL equation with a locally Lipschitz noise because it gives a simple example of SPDE for which
 the difference of two solutions cannot be estimated with the help of only one energy, an essential ingredient
 in
 the proof of Proposition \ref{PropNSFP}.

\subsection{  Notations and Main result }
\

\noindent We set
$$
H=L^2(D;\C),\quad A = -\Delta, \quad  D(A) = H^1_0(D;\C) \cap H^2(D;\C),
$$
and we denote by  $\abs\cdot$, $\abs{\cdot}_p$, $\norm{\cdot}$ and $\norm{\cdot}_s$ the norm of
 $\C$, $L^p(D;\C)$, $H^1_0(D;\C)$ and $H^s(D;\C)$. The norm of $H$ will be denoted by $\abs\cdot$
 when no confusion is possible or $\abs\cdot_H$ otherwise.

\noindent Now we can write problem \eqref{EqIntroCGL} in the form
\begin{equation}\label{EqCGL}
\left\{
\begin{array}{rcl}
du+(\eps+\i)A u \,dt+(\eta+\lambda\i)\abs{u}^{2\sig} u \,dt &=& \phi(u)dW,\\
        u(0)&=& u_0,
\end{array}
\right.
\end{equation}
where $W$ is a cylindrical Wiener process on a Hilbert $U$.

\noindent  In order to have existence and uniqueness of the solution of \eqref{EqCGL},
we make the following assumption

\noindent {\bf H0'} { \it The function $\phi: \,H\to \mathcal{L}_2(U;H)$ is
bounded and local Lypschitz. More precisely, we assume there exists $L>0$
 such that for any $(u_1,u_2)\in H^2$}
$$
\norm{\phi(u_2)-\phi(u_1)}^2_{\mathcal{L}_2(U;H)}\leq L\abs{u_2-u_1}^2
\left(1+\abs{u_1}^{2\sig}+\abs{u_2}^{2\sig}\right).
$$
We set
$$
B_1=\sup_{u\in H}\norm{\phi(u)}^2_{\mathcal{L}_2(U;H)}.
$$
\noindent Under {\bf H0'}, there exists a unique $H$--valued
solution $u(\cdot,W,u_0)$ of \eqref{EqCGL}. This is a continuous
homogenous weak Markov process and a non anticipative measurable
map of $W$. Moreover the law $\Dr(u(\cdot,W,u_0))$ measurably
depends of $u_0$. Existence, uniqueness and measurable dependance
can be established by a contracting fix-point argument applied to
the mild form of \eqref{EqCGL}.

 We denote by $(\mathcal P_t)_{t\in\R^+}$ the Markov
transition semi-group associated to the Markov
 family $\left(u(\cdot,W,u_0)\right)_{u_0\in H}$.

\noindent In our computations, we use the following energy
$$
E_u(t)=\abs{u(t)}^2+\eps \int_{0}^t\norm{u(s)}^2 ds +\eta
\int_{0}^t\abs{u(s)}_{2\sig+2}^{2\sig+2} ds.
$$

\noindent It is well-known  that $\left(A,\Dr(A)\right)$ is
 a self-adjoint operator with discrete spectrum.
We consider $(e_n)_n$ an eigenbasis of $H$ associated to the increasing sequence $(\mu_n)_n$ of eigenvalues
 of $\left(A,\Dr(A)\right)$. We denote by $P_N$ and $Q_N$
the orthogonal projection in  $H$  onto the space $Sp(e_k )_{1\leq N}$ and onto its complementary,
 respectively.

\noindent Now, we state the assumption  which gives the exponential mixing of
$(\Pcal_t)_{t\in \R^+}$ provided it holds for $N$ sufficiently high.

\noindent {\bf H1} {\it  There exists a bounded measurable map $g:H\to \mathcal L(H; U)$ such that for any $u\in H$ }
$$
\phi(u)g(u)=P_N.
$$
\noindent The aim of this section is to establish the following result.
\begin{Theorem}\label{ThCGL}
Assume {\bf H0'}. There exists $N_0(B_1,\eps,\eta,\sig,D,L)$ such that, if  {\bf H1} holds with $N\geq N_0$, then
there exists a unique stationary
 probability measure $\mu$ of $(\Pcal_t)_{t\in \R^+}$ on $H$.
 Moreover, $\mu$ satisfies
\begin{equation}\label{EqThCGLa}
\int_{H} \abs{u}^{2} d\mu(u) < \infty,
\end{equation}
and there exist $C,\gamma>0$  such that for any $\lambda \in \Pcal(H)$
\begin{equation}\label{EqThCGLb}
\|\Pcal^*_t\lambda-\mu\|_*\leq C e^{-\gamma t}\left(1+ \int_{H} \abs{u}^2 d\lambda(u) \right).
\end{equation}
\end{Theorem}
We now prove prove Theorem \ref{ThCGL} by applying Theorem
\ref{ThCrt}.

\subsection{  Energy estimate and Lyapunov structure }
\

Setting $\Hcal=\abs\cdot^2$, the following result states that {\bf
A1} is true.
\begin{Lemma}[The Lyapunov structure]\label{CGLlemNSLyap}
Assume {\bf H0'}. There exist $C_1$ and a family
$(C_\alpha')_\alpha$ only depending on $\eps$, $\eta$, $\sig$,
$B_1$ and $D$ such that
$$
\Espace\left\{\begin{array}{rcl}
\E\left(\abs{u(t,W,u_0)}^{2}\right)&\leq&  e^{-\eps \mu_1  t} \abs{u_0}^{2}+C_1,\\
\E\left(e^{-\alpha
\tau}\abs{u(\tau,W,u_0)}^{2}1_{\tau<\infty}\right)&\leq&
\abs{u_0}^{2}+C_\alpha',
\end{array}\right.$$
for any $t\geq 0$, any $\alpha>0$ and any stopping time $\tau$.
\end{Lemma}

{\it Proof of Lemma  \ref{CGLlemNSLyap}}

\noindent We set $u=u(\cdot,s,W,u_0)$. We apply Ito Formula to
$\abs{u}^2$
$$
d\abs{u}^2+2\eps\norm{u}^2dt+2\eta\abs{u}_{2\sig+2}^{2\sig+2}dt=dM'+\norm{\phi(u)}^2_{\mathcal
L_2(U;H)}dt,
$$
where
 $$
  dM'=2(u,\phi(u)dW).
$$
Hence, taking into account {\bf H0'}, we deduce
\begin{equation}\label{CGLEqA0}
d\abs{u}^2+2\eps\norm{u}^2dt+2\eta\abs{u}_{2\sig+2}^{2\sig+2}dt=dM'+B_1dt.
\end{equation}
Now, we are able to deduce Lemma \ref{CGLlemNSLyap} from
\eqref{CGLEqA0} by applying the same argument we used to deduce
Lemma \ref{lemNSLyap} from \eqref{EqA0}.
 \carre

The next result will be useful in the following.
 \begin{Lemma}[Exponential estimate for the growth of the solution]\label{CGLlemNSnrj}
Assume that {\bf H0} holds. There exists $\gamma_2>0$ only
depending on $\eps$, $\eta$, $\sig$, $B_1$ and $D$ such that
$$
\E\left(\exp\left(\gamma_2\sup_{t\geq
0}\left(E_{u(\cdot,W,u_0)}(t)-B_1 t\right)\right)\right) \leq
2e^{\gamma_2 \abs{u_0}^2},
$$
for any $u_0\in H$.
\end{Lemma}

{\it Proof of Lemma  \ref{CGLlemNSnrj}}

\noindent We set for any $\gamma>0$
$$
\mathcal M_\gamma'(t)=M'(t)-\frac{\gamma}{2}\left<M'\right>(t).
$$
Remarking that
\begin{equation}\label{CGLEqA0.1}
d\left<M'\right>=4\abs{\phi(u)^*u}^2dt\leq 4cB_1\norm{u}^2dt,
\end{equation}
and setting $\gamma_1=\frac{\nu}{2cB_1}$, we deduce from
\eqref{CGLEqA0} that
\begin{equation}\label{CGLEqA2}
E_u(t)\leq \mathcal M_{\gamma_1}(t)+\abs{u_0}^2+ B_1 t.
\end{equation}
Now, we are able to deduce Lemma \ref{CGLlemNSnrj} from
\eqref{CGLEqA2} by applying the same argument we used to deduce
Lemma \ref{lemNSnrj} from \eqref{EqA2}. \carre

\subsection{  Construction of the auxiliary process }\label{Sec3.2}
\

\noindent Now, we build the auxiliary process $\widetilde u$ such
that assumptions {\bf A2},...,{\bf A5} are true. This will allow
to deduce Theorem \ref{ThCGL} from Theorem \ref{ThCrt}.

\noindent Let $K$ be a positive number. We set
$$
F(u)=(\eps+\i) A u+(\eta+\lambda\i)\abs{u}^{2\sig}u,
$$
and
\begin{equation}\label{EqCGL3}
\Espace
\begin{array}{r}
\delta(u,\widetilde u)
=P_N\left((\eta+\lambda\i)\left(\abs{u}^{2\sig}u-\abs{\widetilde
u}^{2\sig}\widetilde u\right)
+ K\abs{P_N\widetilde u}^{2\sig}P_N\left(\widetilde u-u\right)\right)\\
+ L\left(1+\abs{\widetilde
u}^{2\sig}_{H}+\abs{u}^{2\sig}_{H}\right)P_N\left(\widetilde
u-u\right).
\end{array}
\end{equation}
We now consider the following equation
\begin{equation}\label{EqCGLbis}
\Espace
\left\{
\begin{array}{rcl}
d\widetilde u+F(\widetilde u) dt +\delta(u(t,W,u_0),\widetilde u)dt&=& \phi(\widetilde u)dW,\\
                             \widetilde u(0)&=& \widetilde u_0.
\end{array}
\right.
\end{equation}

\noindent It is not difficult to deduce from {\bf H0}  that there
exists a unique $H$--valued solution to \eqref{EqNSbis} when
$(u_0,\widetilde u_0)\in  H^2$. Moreover, $\widetilde u$  is a non
anticipative measurable map of $W$ (See Remark \ref{Rqmes}). It
follows from the uniqueness that $(u,\widetilde u)$ is a
homogenous weak Markov process. Moreover, its law $\Dr(u,
\widetilde u)$ is measurably depending of its initial condition
$(u_0,\widetilde u_0)$.

 We
will denote the dependance of $(u,\widetilde u)$ with respect to
$(t,W,u_0,\widetilde u_0)$ as follows
$$
(u(t),\widetilde u(t))=(u(t,W,u_0),\widetilde u(t,W,u_0,\widetilde
u_0)),
$$
and we deduce {\bf A0}.

\noindent Taking into account {\bf H1} and the uniqueness of the solution of \eqref{EqCGL} under {\bf H0'},
 we deduce {\bf A3} by setting
\begin{equation}\label{EqCGL2}
h(u,\widetilde u)=-g(\widetilde u)\delta(u, \widetilde u).
\end{equation}

We first state a Lemma that will be useful in the following.

\begin{Lemma}[Exponential estimate for the growth of the auxiliary process]\label{lemCGLnrjbis}
Assume {\bf H0'}. There exists $C$, $B$, $K_L$, $\gamma_0'$,
$\gamma_0>0$ not depending on $N$ such that if $K= K_L$
$$
\E\left(\exp\left(\gamma_0'\sup_{t\geq
0}\left(E_{u(\cdot,W,u_0)}(t)+E_{ \widetilde u(\cdot,W,\widetilde
u_0,u_0)}(t)-B t\right)\right)\right) \leq Ce^{\gamma_0
(\abs{u_0}^2+\abs{\widetilde u_0}^2)},
$$
for any $(u_0,\widetilde u_0)\in H^2$.
\end{Lemma}

{\it Proof of Lemma \ref{lemCGLnrjbis}}

\noindent For any function $f$, we denote by $\delta f(u)$ the
value $f(\widetilde u)-f(u)$. Moreover,
 we set $r=\widetilde u-u$.

\noindent Taking the Ito Formula of $\abs{\widetilde u}^2$, we
obtain
\begin{equation}\label{EqA7}
\Espace
\begin{array}{r}
d\abs{\widetilde u}^2 + 2\eps\norm{\widetilde
u}^2dt+2\eta\abs{\widetilde u}^{2\sig+2}_{2\sig+2}dt
+2K\abs{P_N\widetilde u}^{2\sig+2}_{2\sig+2}dt
= 2\left(\widetilde u, \phi(\widetilde u)dW\right)\\
+\norm{\phi(u)}_{\mathcal L_2\left(U;H\right)}^2dt+ I(t)dt
+2\left(P_N \widetilde u,(\eta+\lambda\i)\delta\left(\abs{
u}^{2\sig} u\right)\right)dt,
\end{array}
\end{equation}
where $I=I_1+I_2$ and
$$
\Espace
\begin{array}{lcl}
I_1(t)&=& -2L\left(1+\abs{\widetilde u}^{2\sig}+\abs{u}^{2\sig}\right)\left(P_N \widetilde u,P_N r\right),\\
I_2(t)&=&2K\left(\abs{P_N \widetilde u}^{2\sig}P_N \widetilde
u,P_N u\right).
\end{array}
$$
Applying arithmetico-geometric inequalities, we obtain
\begin{equation}\label{EqA8}
I_2(t)\leq K\abs{P_N\widetilde u}^{2\sig+2}_{2\sig+2}+ c K\abs{
u}^{2\sig+2}_{2\sig+2}.
\end{equation}
Remarking that
$$
\left(P_N \widetilde u,P_N r\right)=\abs{P_N \widetilde
u}^2-\left(P_N \widetilde u,P_N u\right),
$$
we deduce from Schwarz inequality
$$
I_1(t)\leq 2L\left(1+\abs{\widetilde
u}^{2\sig}+\abs{u}^{2\sig}\right)\abs{P_N \widetilde u}\abs{P_N
u},
$$
which implies by applying an arithmetico-geometric inequality
\begin{equation}\label{EqA9}
I_1(t)\leq 1+\frac{K_L}{2}\abs{P_N \widetilde
u}_{2\sig+2}^{2\sig+2}
+\frac{\eta}{2}\left(\abs{u}^{2\sig+2}_{2\sig+2}+\abs{\widetilde
u}^{2\sig+2}_{2\sig+2}\right).
\end{equation}
Applying arithmetico-geometric inequality, it follows
\begin{equation}\label{EqA10}
2\left(P_N \widetilde u,(\eta+\lambda\i)\delta\left(\abs{
u}^{2\sig} u\right)\right)\leq \frac{K_L}{2}\abs{P_N\widetilde
u}^{2\sig+2}_{2\sig+2}
+\frac{\eta}{2}\left(\abs{u}^{2\sig+2}_{2\sig+2}+\abs{\widetilde
u}^{2\sig+2}_{2\sig+2}\right).
\end{equation}
Combining \eqref{CGLEqA0}, \eqref{CGLEqA0.1}, \eqref{EqA7},
\eqref{EqA8}, \eqref{EqA9}, \eqref{EqA10} and {\bf H0'}, we obtain
that if $K=K_L$
\begin{equation}\label{EqA11}
dE_{\widetilde u}+4c_1(1+2K)dE_u\leq\mathcal M(t)+Bt,
\end{equation}
where $\mathcal M$ has been defined as in Section \ref{SecNSnrj}
$$
M(t)=4c_1(1+2K)M'(t)+2\int_s^t\left(u(r),\phi(u(r))dW(r)\right),
\quad \mathcal M(t)=M(t)-\frac{\gamma_0}{2}\left<M\right>(t).
$$
Now, since there is no loss of generality of assuming $c_1\geq 1$,
we are able to deduce Lemma \ref{lemCGLnrjbis} from \eqref{EqA11}
by applying the same argument we used to deduce Lemma
\ref{lemNSnrj} from \eqref{EqA2}. \carre

We now fix $K=K_L$ and we state a result analogous to Proposition
\ref{PropNSFP}.
\begin{Proposition}\label{PropCGLFP}
 Assume that {\bf H0'} holds. There exist $(C_N)_N$, $\alpha$, $\gamma>0$,
$\gamma'\in (0,1]$ and $N_0$ depending only on $B_1$, $L$, $\eps$,
$\eta$ and $D$ such that
 for any  $N\geq N_0$ and any $(t,u_0,\widetilde u_0)\in (0,\infty)\times H^2$
$$
\E\left(\left(e^{ \frac{\eps\mu_1}8 t}\abs{r(t)}^2+\int_0^t
e^{\frac{\eps\mu_1}8s}Z(s)ds\right)^{\gamma'}\right) \leq
C_{N}\abs{r(0)}^{2\gamma'}e^{\gamma (\abs{u_0}^2+\abs{\widetilde
u_0}^2)},
$$
where
$$
\Espace\left\{\begin{array}{rcl}
 r&=&\widetilde u-u,\quad \widetilde
u=\widetilde u(\cdot,W,u_0,\widetilde u_0),\quad
u=u(\cdot,W,u_0),\quad\beta= 4\sig+\frac{2\sig+2}{2-\sig}\\
Z&=&(\norm{u}^2+\norm{\widetilde
u}^2+\abs{u}_{2\sig+2}^{2\sig+2}+\abs{\widetilde
u}_{2\sig+2}^{2\sig+2})\abs{r}^2+(1+\abs{u}^\beta+\abs{\widetilde
u}^\beta)\norm{r}^2.
\end{array}\right.$$
\end{Proposition}
By Chebyshev inequality,  {\bf A2} immediately follows.

{\it Proof of Proposition \ref{PropCGLFP}}

\noindent For any function $f$, we denote by $\delta f(u)$ the
value $f(\widetilde u)-f(u)$.
 Taking the difference between \eqref{EqCGL} and \eqref{EqCGLbis}, we obtain
$$
\Espace
\begin{array}{r}
dr + (\eps +\i)A r dt + L  \left(1+\abs{\widetilde
u}^{2\sig}_{H}+\abs{u}^{2\sig}_{H}\right)P_N r dt
+KP_N\left(\abs{P_N\widetilde u}^{2\sig}P_Nr\right)dt
\quad\quad\\
=\delta \phi(u)dW
-(\eta+\lambda\i)Q_N\delta\left(\abs{u}^{2\sig}u\right)dt.
\end{array}
$$
Hence, applying Ito Formula to $\abs{r}^2$, we have
\begin{equation}\label{EqA4}
\Espace
\begin{array}{r}
d\abs{r}^2 + 2\eps\norm{r}^2dt +2L  \left(1+\abs{\widetilde
u}^{2\sig}+\abs{u}^{2\sig}\right)\abs{P_N r}^2dt
\leq 2\left(r,\delta \phi(u)dW\right)\\
-2\left(Q_Nr,(\eta+\lambda\i)\delta\left(\abs{u}^{2\sig}u\right)\right)dt
+\norm{\delta\phi(u)}_{\mathcal L_2\left(U;H\right)}^2dt.
\end{array}
\end{equation}
Remarking that for any $(x,y)\in \C^2$
\begin{equation}\label{EqAC}
\abs{\abs{x}^{2\sig}x-\abs{y}^{2\sig}y}\leq
C_\sig\left(\abs{x}^{2\sig}+\abs{y}^{2\sig}\right)\abs{x-y}^2,
\end{equation}
it follows from H\"older inequality
$$
-\left(Q_Nr,(\eta+\lambda\i)\delta\left(\abs{u}^{2\sig}u\right)\right)
\leq c\abs{Q_Nr}_{2\sig+2}\abs{r}_{2\sig+2}\left(\abs{\widetilde
u}_{2\sig+2}^{2\sig}+\abs{u}_{2\sig+2}^{2\sig}\right).
$$
Setting
$$
s_0=\frac{\sig d}{2\sig+2}, \quad s_+=\frac{1}{\sig+1},
$$
we deduce from Sobolev Embedding $H^{s_0}\left(D;\C\right)\subset
L^{2\sig+2}\left(D;\C\right)$ that
$$
-\left(Q_Nr,(\eta+\lambda\i)\delta\left(\abs{u}^{2\sig}u\right)\right)
\leq c\norm{Q_Nr}_{s_0}\norm{r}_{s_0}\left(\abs{\widetilde
u}_{2\sig+2}^{2\sig}+\abs{u}_{2\sig+2}^{2\sig}\right).
$$
Remarking that $s_0<s_+<1$ and that
$\norm{Q_Nr}_{s_+}\leq\mu_{N+1}^{\frac{s_+-s_0}{2}}\norm{Q_Nr}_{s_0}$,
we
 deduce from interpolatory inequality
$$
-\left(Q_Nr,(\eta+\lambda\i)\delta\left(\abs{u}^{2\sig}u\right)\right)
\leq
c\mu_{N+1}^{-\frac{s_+-s_0}{2}}\norm{r}^{2s_+}\abs{r}^{2(1-s_+)}
\left(\abs{\widetilde
u}_{2\sig+2}^{2\sig}+\abs{u}_{2\sig+2}^{2\sig}\right).
$$
Hence, it follows from arithmetico-geometric inequality that there
exists $\alpha\in (0,1)$
 only depending on $\sig$ and $d$ such that
\begin{equation}\label{EqA5}
-\left(Q_Nr,(\eta+\lambda\i)\delta\left(\abs{u}^{2\sig}u\right)\right)
\leq
\frac{\eps}{2}\norm{r}^{2}+\frac{c}{\mu_{N+1}^\alpha}\abs{r}^{2}
\left(\abs{\widetilde
u}_{2\sig+2}^{2\sig+2}+\abs{u}_{2\sig+2}^{2\sig+2}\right).
\end{equation}
Recall {\bf H0'},
$$
\norm{\delta\phi(u)}_{\mathcal L_2\left(U;H\right)}^2\leq
L\left(1+\abs{\widetilde
u}^{2\sig}+\abs{u}^{2\sig}\right)\abs{r}^2.
$$
Remarking that
$$
\abs{r}^2\leq \abs{P_Nr}^2+\frac{1}{\mu_{N+1}^{s_+}}\norm{Q_N
r}^2_{s_+},
$$
 and making interpolatory inequality analogous to those done to obtain \eqref{EqA5}, we obtain
\begin{equation}\label{EqA6}
\Espace
\begin{array}{r}
\norm{\delta\phi(u)}_{\mathcal L_2\left(U;H\right)}^2\leq
\frac{\eps}{4}\norm{r}^{2}+\frac{c}{\mu_{N+1}^\alpha}\abs{r}^{2}
\left(\abs{\widetilde
u}_{2\sig+2}^{2\sig+2}+\abs{u}_{2\sig+2}^{2\sig+2}\right)
\quad\quad\quad\quad\quad\quad\\
+L\left(1+\abs{\widetilde
u}^{2\sig}+\abs{u}^{2\sig}\right)\abs{P_Nr}^2.
\end{array}
\end{equation}
Hence, setting
$$
V(t)=1+\abs{\widetilde u(t)}_{2\sig+2}^{2\sig+2}
+\abs{u(t)}_{2\sig+2}^{2\sig+2},
$$
we deduce from \eqref{EqA4}, \eqref{EqA5} and \eqref{EqA6} that
there exists $\Lambda_1>0$ such that
$$
\Espace
\begin{array}{r}
d\abs{r}^2 + \left( \frac{5}{4}\eps\norm{r}^2
-\frac{\Lambda_1}{\mu_{N+1}^\alpha}V \abs{r}^2
+L\left(1+\abs{\widetilde u}^{2\sig}+\abs{u}^{2\sig}\right)\abs{P_Nr}^2\right)dt\quad\quad\quad\quad\\
\leq 2\left(r,\delta \phi(u)dW\right).
\end{array}
$$
Taking into account \eqref{EqA11}, Ito Formula of
$\left(\mu_{N+1}^\alpha+4c_1(1+2K)\abs{u}^2+\abs{\widetilde u}^2
\right) \abs{r}^2$ gives
$$
\Espace
\begin{array}{r}
d(X\abs{r}^2) +
\left(\frac{\eps\mu_1}4-\frac{\Lambda_1}{\mu_{N+1}^\alpha}V
\right)X\abs{r}^2dt+ \eps X\norm{r}^2dt+Y\abs{r}^2dt
\quad\quad\quad\quad\\
\leq dM_1+B\abs{r}^2dt+d\left<M,M_r\right>.
\end{array}
$$
where
$$
\Espace\left\{
\begin{array}{rcl}
M_r(t)&=&\int_0^t\left(r,\delta\phi(u)dW\right),\\
X(t)&=&\mu_{N+1}^\alpha+4c_1(1+2K)\abs{u(t)}^2+\abs{\widetilde u(t)}^2,\\
Y(t)&=&\mu_{N+1}^\alpha+4c_1(1+2K)\left(\eps\norm{u(t)}^2+\eta\abs{u(t)}_{2\sig+2}^{2\sig+2}\right)+
\left(\eps\norm{\widetilde u(t)}^2+\eta\abs{\widetilde
u(t)}_{2\sig+2}^{2\sig+2}\right),\\
dM_1&=&XdM_r+\abs{r}^2dM.
\end{array}\right.
$$
Notice that
$$
d\left<M,M_r\right>\leq
c\left(\abs{\phi^*(u)u}+\abs{\phi^*(\widetilde u)\widetilde
u}\right)\abs{\delta\phi^*r}dt
$$
It follows from {\bf H0'} that
$$
d\left<M,M_r\right>\leq c\left(\abs{u}+\abs{\widetilde u}\right)
\abs{r}^2\left(1+\abs{u}^\sig+\abs{\widetilde u}^\sig\right),
$$
which yields
$$
B\abs{r}^2dt+d\left<M,M_r\right>\leq \frac{c}{\mu_{N+1}^\alpha}
V\abs{r}^2dt\leq \frac{c}{\mu_{N+1}^\alpha} VX\abs{r}^2dt.
$$
It follows that
$$
d(X\abs{r}^2) +
\left(\frac{\eps\mu_1}4-\frac{\Lambda_2}{\mu_{N+1}^\alpha}V
\right)(X\abs{r}^2)dt+ \eps X\norm{r}^2dt+Y\abs{r}^2dt \leq dM_1.
$$
Integrating this formula and taking the expectation, it follows
$$
\Espace \E\left(
\begin{array}{l}
e^{\frac{\eps\mu_1}8t}G(t)^{-1}X(t)\abs{r(t)}^2+\\
\int_0^t e^{\frac{\eps\mu_1}8s}G(s)^{-1}\left(\eps
X(s)\norm{r(s)}^2+Y(s)e^{\frac{\Lambda}{\mu_{N+1}^\alpha}
(\abs{u(s)}^2+\abs{\widetilde u(s)}^2)}\abs{r(s)} \right)ds
\end{array}
\right)
\leq \abs{r(0)}^{2},
$$
where
$$
G(t)=e^{
-\frac{\eps\mu_1}8t+\frac{\Lambda}{\mu_{N+1}^\alpha}(E_u(t)+E_{\widetilde
u}(t))}.
$$
Notice that  $1+\abs{u}^{\beta}+\abs{\widetilde u}^\beta\leq
C_{N}Y e^{\frac{\Lambda}{\mu_{N+1}^\alpha}(
\abs{u}^2+\abs{\widetilde u}^2)}$.
 Hence, choosing $N$ sufficiently high and $\eps>0$ sufficiently small, we
are able to deduce Proposition \ref{PropCGLFP} from Lemma
\ref{lemCGLnrjbis} as we have deduced Proposition \ref{PropNSFP}
from Lemma \ref{lemNSnrj}.
 \carre

 To end the proof of Theorem \ref{ThCGL}, we establish that
 \begin{equation}\label{Fin}
\abs{h(u,\widetilde u)}^2\leq C_N Z.
\end{equation}
Hence, by Chebyshev inequality,  we deduce {\bf A4} and {\bf A5}
 from Proposition \ref{PropCGLFP}.

{\it Proof of \eqref{Fin}}

Taking into account {\bf H1} and \eqref{EqCGL2},
 we remark that it is sufficient to establish
$$
 \abs{\delta(\widetilde u,u)}^2\leq C_N Z.
$$
Recalling   \eqref{EqCGL3}, we deduce that it is sufficient to estimate the three following values
$$
\Espace
\begin{array}{lcl}
I_1&=& \abs{P_N\left(\abs{u}^{2\sig}u-\abs{\widetilde u}^{2\sig}\widetilde u\right)}^2,\\
I_2&=&\abs{P_N\left( \abs{P_N\widetilde u}^{2\sig}P_N r\right)}^2,\\
I_3&=& \left(1+\abs{\widetilde
u}^{4\sig}+\abs{u}^{4\sig}\right)\abs{P_N r}^2.
\end{array}
$$
For $I_3$, the result is obvious.
 Applying successively H\"older inequality and the equivalence of the norm in finite-dimensional spaces, it follows
$$
 I_2=\abs{P_N\left( \abs{P_N\widetilde u}^{2\sig}P_N r\right)}^2
\leq  \abs{P_N\widetilde u}^{4\sig}_{4\sig+2}\abs{P_N
r}^{2}_{4\sig+2}\leq K_N'\abs{\widetilde u}^{4\sig}\abs{r}^2\leq
C_N Z.
$$
The equivalence of the norm in finite-dimensional spaces gives
$$
I_1=\abs{P_N\left(\abs{u}^{2\sig}u-\abs{\widetilde
u}^{2\sig}\widetilde u\right)}^2 \leq
K_N'\abs{\abs{u}^{2\sig}u-\abs{\widetilde u}^{2\sig}\widetilde
u}_1^2.
$$
Hence, we deduce from \eqref{EqAC} that
$$
I_1\leq K_N'\abs{\left(\abs{u}^{2\sig}+\abs{\widetilde
u}^{2\sig}\right)r}_1^2.
$$
We first treat the case $\sig\leq 1$. In that case, a Schwarz
inequality gives
$$
I_1 \leq K_N'\left(\abs{u}^{4\sig}_{4\sig}+\abs{\widetilde
u}^{4\sig}_{4\sig}\right)\abs{r}^2 \leq C_N Z.
$$
Now, we treat the case $\sig\in (1,2)$. In that case $d=1$.
Sobolev Embedding $H^1(D)\subset L^\infty(D)$ gives
$$
I_1\leq K_N'\left(\abs{u}_{2\sig}^{4\sig}+\abs{\widetilde
u}_{2\sig}^{4\sig}\right)\norm{r}^2.
$$
Sobolev Embedding $H^{\frac{\sig-1}{2\sig}}(D)\subset
L^{2\sig}(D)$ and interpolatory inequality gives
$$
\abs{u}_{2\sig}^{4\sig}\leq
\norm{u}_{\frac{\sig-1}{2\sig}}^{4\sig}\leq
\norm{u}^{2\sig-2}\abs{u}^{2\sig+2}\leq
\norm{u}^2+\abs{u}^{\frac{2\sig+2}{2-\sig}},
$$
which yields
$$
I_1 \leq C_N Z.
$$
So we have $I_i\leq C_N Z$ for $i=1,2,3$, which yields
\eqref{Fin}.
 \carre

\subsection{Conclusion and Remarks}\label{sec_Conclusion_CGL}
 \

We have proved that assumptions of Theorem \ref{ThCrt} are
verified. So Theorem \ref{ThCGL} follows.

Actually, {\bf A0} is consequence of the well-posedness of
equations.

We set $\Hcal=\abs\cdot^2$.  As shown in Lemma \ref{CGLlemNSLyap},
we have {\bf A1}.

Assumption {\bf A3}  has been proved at the beginning of
 section 4.3.

We deduce {\bf A2} directly from Proposition \ref{PropCGLFP} and
{\bf A4}, {\bf A5} from \eqref{Fin} and Proposition
\ref{PropCGLFP}.

  Since {\bf A0},\dots,
{\bf A5} are established, we can apply Theorem \ref{ThCrt} which
yields Theorem \ref{ThCGL}.

\begin{Remark}
 Notice that the proof Proposition \ref{PropCGLFP} is much more
 difficult than the proof of Proposition \ref{PropNSFP} because
  it involves energies control of both $u_1$ and $\widetilde u$.
   Then, in
order to apply it, we have to establish energy estimate of both
$u_1$ and $\widetilde u$ (Lemma \ref{lemCGLnrjbis}). That is the
reason why
 the building of the auxiliary process is more complicated than
  for NS or for CGL with a global Lipschitz noise.
  For instance, if we set
$$
\delta(\widetilde u, u_1)= KP_N(\widetilde u-u_1),
$$
then we would have an energy Lemma on $\widetilde u$ analogous to
Lemma \ref{lemCGLnrjbis} and
 a Proposition analogous to Proposition \ref{PropCGLFP} under a truncation
 condition on $\widetilde u$. The
 problem is that we do not know if it is possible to choose
$K$ and $N$ such that both result are true. So we can not combine
them.
\end{Remark}

\begin{Remark}\label{RqNew2}
There exists a lot of variations of Theorem \ref{ThCGL}. Hence
\eqref{EqThCGLa} could be strengthened into
$$
\int_{H}
\exp\left(\alpha_1(B_1,\eps,\eta,D)\abs{u}^2\right)d\mu(u) <
\infty.
$$
\noindent Moreover, one can work in $H^1_0(D;\C)$ in the
defocusing case under the $H^1$--subcritical condition
$(d-2)\sig<2$ and in the focusing
 case under the subcritical condition $\sig<\frac{2}{d}$. In the focusing case, the rate of convergence is
greater than any power of time instead of being exponential in
time.

\noindent In {\bf H1}, the boundedness of $g$ could be replaced by
the existence of $C$ such that
$$
\abs{g(u)}_{\mathcal L(H;U)}\leq
C\exp\left(\frac{\eps\mu_1}{4}\abs{u}^2\right).
$$
Contrary to Navier--Stokes,  the coefficient in the exponential
cannot be as high as we want because it seems that for the locally
Lipschitz CGL, there is no property analogous to Proposition
\ref{PropNSFP}.

\noindent In {\bf H0'}, the boundedness of $\phi$ could be
replaced by $\abs{\phi(u)}\leq C(1+\abs{u}^\gamma)$. If
$\gamma\leq \sig$, then the rate of convergence remains
exponential. If $\gamma<\sig+1$, then the rate of convergence
becomes greater than any power of time instead of being
exponential in time. Moreover for any $p$ there exists $c_p$ such
that if there exists $C$ such that $\abs{\phi(u)}\leq
C+c_p\abs{u}^{\sig+1}$, then the rate of convergence is greater
than $(1+t)^{-p}$.
 \end{Remark}





\footnotesize

\begin{thebibliography}{9}




\BLANC{
\bibitem{Marc} M. Barton-Smith,
{\it Etudes th\'eoriques et num\'eriques sur les \'equations de Schr\"odinger
 et Ginzburg-Landau stochastiques },
 Th\`ese soutenue \`a Orsay, 2002.
}

\bibitem{Marc} M. Barton-Smith,
{\it Invariant measure for the stochastic Ginzburg Landau equation},
 Nonlinear Differential Equations Appl. {\bf 11}  ,  no. 1, 29--52, 2004.

\bibitem{Bebouche} P. Bebouche, A. J\"ungel,
{\it Inviscid limits of the Complex Ginzburg--Landau Equation},
Commun. Math. Phys. {\bf 214}, 201-226, 2000.




\bibitem{BKL} J. Bricmont, A. Kupiainen and R. Lefevere,
{\it Exponential mixing for the 2D stochastic Navier-Stokes dynamics},
Commun. Math. Phys. {\bf 230}, No.1, 87-132, 2002.

\bibitem{Cons} P. Constantin and C. Foias,
{\it Navier-Stokes Equations},
University of Chicago Press, Chicago, IL, 1988.

\bibitem{DebusscheNS3D}  G. Da Prato, A. Debussche,
{\it Ergodicity for the 3D stochastic Navier-Stokes equations}
J. Math. Pures Appl. (9)  82  (2003),  no. 8, 877--947.


\bibitem{DPZ1} G. Da Prato and J. Zabczyk, {\it Stochastic equations
in infinite
dimensions,}
   Encyclopedia of Mathematics and its Applications,  Cambridge University
Press, 1992.


\bibitem{DPZ2} G. Da Prato and J. Zabczyk,  {\it Ergodicity for Infinite Dimensional Systems,} London
Mathematical Society Lecture Notes, n.229, Cambridge University Press, 1996.



\bibitem{ODASSO2} A. Debussche, C. Odasso,
{\it Ergodicity for the weakly damped stochastic Non-linear
Schr\"odinger  equations}, to appear in Journal of Evolution
Equations.


\bibitem{EMS} W. E, J.C. Mattingly, Y. G. Sinai,
{\it Gibbsian dynamics and ergodicity for the stochastically forced
Navier-Stokes equation},
Commun. Math. Phys. {\bf 224}, 83--106, 2001.



\bibitem{FKLT}
{\sc G.E. Falkovich, I. Kolokolov, V. Lebedev, S.K. Turitsyn},
{\em Self-focusing in the perturbed and unperturbed nonlinear
Schroedinger equation in critical dimension},
J. Appl. Math. {\bf 60}, 183--240 (2000).


\bibitem{FlandoliMaslowsky}  F. Flandoli and B. Maslowski,
     {\it Ergodicity of the 2--D
    Navier--Stokes equation
    under random perturbations}, Commun. Math. Phys. {\bf 171}, 119--141, 1995.

\bibitem{GL} V. Ginzburg, L. Landau,
{\it On the theorie of superconductivity},
Zh. Eksp. Fiz. {\bf 20}, 1064(1950) English transl. in: {\it Men of Physics: L.D. Landau.} Vol. I. Ter Haar (ed.).
 New York: Pergammon Press. 1965. pp. 546-568


\bibitem{goubet} O. Goubet, {\it Regularity of the attractor for a weakly damped nonlinear
Schr\"odinger equation}, Appl. Anal. 60, No.1-2, 99-119, 1996.

\bibitem{H} M. Hairer,
{\it Exponential Mixing Properties of Stochastic PDEs Through
  Asymptotic   Coupling},
Proba. Theory Related Fields,{\bf 124}, 3 :345-380, 2002.

\bibitem{HM04} M. Hairer, J. Mattingly, {\it Ergodicity of the 2D Navier-Stokes equations with
degenerate forcing}, preprint.

\bibitem{Huber} G. Huber, P. Alstrom,
{\it Universal Decay of vortex density in two dimensions},
Physica A {\bf 195}, 448-456, 1993.


\bibitem{K}  S. Kuksin,
{\it On exponential convergence to a stationary mesure for nonlinear
 PDEs},
The  M. I. Viishik Moscow PDE seminar, Amer. Math. Soc. Trans. (2), vol 206,
 Amer. Math. Soc., 2002.


\bibitem{KS2000} S. Kuksin and A. Shirikyan,
{\it Stochastic dissipative PDE's and Gibbs measures},
 Comm. Math. Phys.{\bf  213} (2000), no. 2, 291--330.


\bibitem{KS00}  S. Kuksin, A. Shirikyan,
{\it Stochastic dissipative PDE's  and Gibbs measures},
Commun. Math. Phys. {\bf 213}, 291--330, 2000.



\bibitem{KS1}  S. Kuksin and A. Shirikyan,{\it Ergodicity for the randomly forced 2D Navier-Stokes equations},
 Math. Phys. Anal. Geom. {\bf 4}, 2001.

\bibitem{KS}  S. Kuksin, A. Shirikyan,
{\it A coupling approach to randomly forced randomly forced PDE's I},
Commun. Math. Phys. {\bf 221}, 351--366, 2001.


\bibitem{KS2}  S. Kuksin, A. Piatnitski, A. Shirikyan,
{\it A coupling approach to randomly forced randomly forced PDE's II},
Commun. Math. Phys. {\bf 230}, No.1, 81-85, 2002.


\bibitem{KS3}  S. Kuksin, A. Shirikyan,
{\it Coupling approach to white-forced nonlinear
  PDEs},
  J. Math. Pures Appl.  1 (2002) pp. 567-602.



\bibitem{KSlimit}  S. Kuksin, A. Shirikyan,
{\it Randomly forced CGL equation: Stationnary measure and the inviscid limit},
J. Phys. A {\bf 37},  no. 12, 3805--38222004.


\bibitem{Lindvall} T. Lindvall,
{\it Lectures on the coupling method}, John Wiley and Sons, New
York, 1992.



\bibitem{Matt} J. Mattingly,
{\it Exponential convergence for the stochastically forced
  Navier-Stokes equations and other partially dissipative dynamics},
Commun. Math. Phys. {\bf 230}, 421-462, 2002.

\bibitem{MattPar}  J. Mattingly, E. Pardoux,
{\it Ergodicity of the 2D Navier-Stokes Equations with Degenerate
Stochastic Forcing },
 preprint 2004.





\bibitem{Newel1} A. Newel, J. Whitehead,
{\it Finite bandwidth, finite amplitude convection},
J. Fluid Mech. {\bf 38}, 279-303, 1969.

\bibitem{Newel2} A. Newel, J. Whitehead,
{\it Review of the finite bandwidth concept},
H. Leipholz. editor. Proceedings of the Internat. Union of Theor. and Appl. Math.. Berlin: Springer. 1971, pp 284-289.

\bibitem{ODASSO1} C. Odasso,
{\it Ergodicity for the stochastic Complex Ginzburg--Landau
equations}, to appear in Annales de l'institut Henri-Poincaré,
Probabilités et Statistiques.


\bibitem{S}   A. Shirikyan,
{\it Exponential mixing for 2D Navier-Stokes equation pertubed by an unbounded
 noise},
J. Math. Fluid Mech. {\bf 6},  no. 2, 169--193, 2004.

\bibitem{Temam}   R. Temam,
{\it Navier--Stokes Equations. Theory and Numerical Analysis.},
North-Holland, Amsterdam-New York-Oxford,1977.


\end{thebibliography}
\end{document}